\documentclass[12pt]{article}
\usepackage{amsmath,amsfonts}
\usepackage{verbatim}
\usepackage{latexsym}
\usepackage{graphicx}
 \usepackage{color}
\makeatletter \@addtoreset{equation}{section}

\begin{document}

\newcommand{\E}{\mathbb{E}}
\newcommand{\PP}{\mathbb{P}}
\newcommand{\RR}{\mathbb{R}}

\newcommand{\Dt}{\D t}
\newcommand{\bX}{\bar X}
\newcommand{\bx}{\bar x}
\newcommand{\by}{\bar y}
\newcommand{\bp}{\bar p}
\newcommand{\bq}{\bar q}
\newcommand{\bi}{\bar i}

\newtheorem{theorem}{Theorem}[section]
\newtheorem{lemma}[theorem]{Lemma}
\newtheorem{coro}[theorem]{Corollary}
\newtheorem{defn}[theorem]{Definition}
\newtheorem{assp}[theorem]{Assumption}
\newtheorem{expl}[theorem]{Example}
\newtheorem{prop}[theorem]{Proposition}
\newtheorem{rmk}[theorem]{Remark}

\newcommand\tq{{\scriptstyle{3\over 4 }\scriptstyle}}
\newcommand\qua{{\scriptstyle{1\over 4 }\scriptstyle}}
\newcommand\hf{{\textstyle{1\over 2 }\displaystyle}}
\newcommand\athird{{\scriptstyle{1\over 3 }\scriptstyle}}
\newcommand\hhf{{\scriptstyle{1\over 2 }\scriptstyle}}

\newcommand{\eproof}{\indent\vrule height6pt width4pt depth1pt\hfil\par\medbreak}

\def\a{\alpha} \def\g{\gamma}
\def\e{\varepsilon} \def\z{\zeta} \def\y{\eta} \def\o{\theta}
\def\vo{\vartheta} \def\k{\kappa} \def\l{\lambda} \def\m{\mu} \def\n{\nu}
\def\x{\xi}  \def\r{\rho} \def\s{\sigma}
\def\p{\phi} \def\f{\varphi}   \def\w{\omega}
\def\q{\surd} \def\i{\bot} \def\h{\forall} \def\j{\emptyset}

\def\be{\beta} \def\de{\delta} \def\up{\upsilon} \def\eq{\equiv}
\def\ve{\vee} \def\we{\wedge}

\def\F{{\cal F}}
\def\T{\tau} \def\G{\Gamma}  \def\D{\Delta} \def\O{\Theta} \def\L{\Lambda}
\def\X{\Xi} \def\S{\Sigma} \def\W{\Omega}
\def\M{\partial} \def\N{\nabla} \def\Ex{\exists} \def\K{\times}
\def\V{\bigvee} \def\U{\bigwedge}

\def\1{\oslash} \def\2{\oplus} \def\3{\otimes} \def\4{\ominus}
\def\5{\circ} \def\6{\odot} \def\7{\backslash} \def\8{\infty}
\def\9{\bigcap} \def\0{\bigcup} \def\+{\pm} \def\-{\mp}
\def\[{\langle} \def\]{\rangle}

\def\tl{\tilde}
\def\trace{\hbox{\rm trace}}
\def\diag{\hbox{\rm diag}}
\def\for{\quad\hbox{for }}
\def\refer{\hangindent=0.3in\hangafter=1}

\newcommand\wD{\widehat{\D}}

\thispagestyle{empty}

\title{
\bf The Partially Truncated Euler--Maruyama Method for Super-linear Stochastic Delay Differential
Equations with Variable Delay and Markovian Switching}

\author{Yuhao Cong\textsuperscript{a}\quad
Weijun Zhan\textsuperscript{a}\thanks{Corresponding author, Email: weijunzhan@hotmail.com}\quad
Qian Guo\textsuperscript{b}\\
\textsuperscript{a}Department of Mathematics,\\ Shanghai University, Shanghai 200444, China\\
\textsuperscript{b}Department of Mathematics,\\ Shanghai Normal University, Shanghai, China, 200234
}
\date{}
\maketitle

\begin{abstract}

A class of super-linear stochastic delay differential equations (SDDEs) with variable delay and Markovian switching is considered.
The main aim of this paper is to investigate the convergence and stability properties of partially truncated Euler-Maruyama (EM)
method applied to the SDDEs with variable delay and Markovian switching under the generalized Khasminskii-type condition.


\medskip\noindent
{\small\bf Key words: } partially truncated Euler--Maruyama, stochastic delay differential equations, variable delay, Markovian switching, super-linear.

\medskip \noindent
{\small\bf  Mathematical Subject Classifications  (2010)}:  60H35, 60J10.
\end{abstract}

\section{Introduction}

Systems in many branches of science and industry do not only depend on the present state but also on the past ones. Therefore, stochastic delay differential equations with constant delay have been extensively used to model such systems \cite{MR05,M2002,M02}. However, it is more realistic for the SDDEs including varying time delay. Forever, these systems may often experience abrupt changes by their structure and parameters and continuous-time Markovian chains. Hence the SDDEs with variable delay and Markovian switching has more meaningful for practice. The sufficient conditions for existence and uniqueness of the exact solution for SDDEs with Markovian were given under the local Lipschitz condition in \cite{Li2010,MY06}. Without the local Lipschitz condition and the linear growth condition, Wang et al. \cite{5} investigated the existence and uniqueness of solutions to the nonlinear SDDEs with Markovian switching. These different conditions show that the coefficients are controlled by polynomial functions.


Generally speaking, explicit solutions or the probability distribution of solutions are not known, the numerical solution for
SDDEs  has received a great deal of attention \cite{BB,E.B2000,MS,M11,Mao2007,wu2010}.
There is much literature concerned with the numerical approximation to SDDEs with Markovian switching. Here we just mention a few of them. On the one hand,
under the global Lipschitz condition, Mao and Yuan \cite{15} proved that EM method
can converge to the exact solution, as well as in \cite{19} presented the EM method can preserve converge to the linear SDDEs with Markovian switching.
Li et al. \cite{Li2010} designed the EM method can converge to the true solutions under the local Lipschitz and Lyapunov-type conditions.
Further, Wang et al. \cite{9} investigated the numerical solutions of SDDEs with Markovian switching by Taylor approximations.
On the other hand, some other methods have been considered to deal with the SDDEs with Markovian switching. For example,
Liu et al. \cite{20} designed semi-implicit Euler method for linear SDDEs. Further, when the drift coefficient obeys the one-sided Lipschitz
and polynomial growth conditions, and the diffusion coefficient is polynomial growth, the convergence in probability of BEM solution was proved in \cite{18}.

Meanwhile, Mao et al. \cite{1,2} showed that the exponential stability
of the exact solution with local Lipschitz and linear growth conditions. For the nonlinear SDDEs with Markovian switching, Mao \cite{7} gave the asymptotic
stability under some specified conditions. Further, Zhu et al. \cite{3} investigated some novel sufficient conditions to  $p$-th moment exponential stability
by using a Lyapunov function and generalized Halanay inequality.

To our knowledge, there is little literature concerned with the explicit method applied to super-linear SDDEs with Markovian switching in the strong sense.
Recently, Mao \cite{mao2015truncated,mao2016convergence} developed a new explicit numerical method, called the truncated EM method, for SDEs under
the Khasminskii-type condition plus the local Lipschitz condition and established the strong convergence theory. Guo et al. \cite{GLMY2017} modified
the truncated EM method by using a partially truncated technique so that the numerical solution can preserve the mean square exponential stability
and asymptotic boundedness of analytical solution to the underlying SDEs. Moreover, Guo et al. \cite{12} and Zhang et al. \cite{22} discussed the convergence
of SDDE with constant delay by truncated EM method and partially truncated EM method under the Khaminskii-type condition, respectively.
Yet, they did not concern with the SDDEs with variable time delay and Markovian switching.
In this paper, we will develop the partially truncated EM method for the super-linear SDDEs with variable time delay and Markovian switching under the generalized
Khasminskii-type condition, meanwhile we give the strong convergence rate and almost sure stability of the numerical solution.

This paper is organized as follows: we will introduce the necessary notion, state the generalized Khasminskii-type condition and give the definition of the partially truncated EM
method for SDDEs with variable delay and Markovian switching in section 2. We will establish the strong convergence theory and convergence rate,
meanwhile, illustrate the theory by examples in sections 3. In section 4, we will present the almost sure exponential stability theory and give some examples.
\section{The Partially Truncated Euler-Maruyama Method }

Throughout this paper, unless otherwise specified, we use the
following notation. Let $|\cdot|$ be the Euclidean norm in $\RR^n$. If
$A$ is a vector or matrix, its transpose is denoted by $A^T$. If $A$
is a matrix, its trace norm is denoted by $|A|=\sqrt{ \trace(A^TA)
}$. Let $\RR_+=[0,\8)$ and $\T >0$.
Let $(\W , {\cal F}, \{{\cal F}_t\}_{t\ge 0}, \PP)$ be a complete
probability space with a filtration  $\{{\cal F}_t\}_{t\ge 0}$
satisfying the usual conditions (i.e., it is increasing and right
continuous while ${\cal F}_0$ contains all $\PP$-null sets). Let
$B(t)= (B_1(t), \cdots, B_m(t))^T$ be an $m$-dimensional Brownian
motion defined on the probability space.
Moreover, for two real numbers $a$ and $b$, we use $a\ve b=\max(a,b)$
and $a\we b=\min(a,b)$. If $G$ is a set, its indicator function is denoted by
$I_G$, namely $I_G(x)=1$ if $x\in G$ and $0$ otherwise.
If $a$ is a real number, we denote by $\lfloor a\rfloor$ the largest integer
which is less or equal to $a$, e.g., $\lfloor -1.2\rfloor = -2$ and $\lfloor 2.3\rfloor =2$.

Let $r(t)$ be a right-continuous Markov chain on the probability space taking values in
a finite state space $\mathbb{S}=\{1,2,\ldots,N \}$ with the generator $\Gamma=(\gamma_{ij})_{N\times N}$
gives by
\begin{equation*}
P\{r(t+\D)=j|r(t)=i\}=
\begin{cases}
\gamma_{ij}+o(\D) &\quad \text{if} \quad i\neq j\\
1+\gamma_{ij}+o(\D)& \quad \text{if} \quad i=j
\end{cases}
\end{equation*}
where $\D>0$. Here $\gamma_{ij}\geq 0$ is the transition rate from $i$ to $j$ if $i\neq j$ while
$$
\gamma_{ii}=-\sum_{i\neq j}\gamma_{ij}.
$$
We assume that the Markov chain $r(\cdot)$ is independent of the Brownian motion $B(\cdot)$. It is well known
that almost every sample path of $r(\cdot)$ is a right-continuous step function with finite number of simple jumps
in any finite subinterval of $\mathbb{R}_{+}$.

Consider a nonlinear stochastic differential equation with variable delay and Markovian switching of the form
\begin{eqnarray}\label{SDDE}
dx(t)& =& f(x(t),x(t-\delta(t)),r(t))dt+g(x(t), x(t-\delta(t)),r(t))dB(t), \quad t\ge 0,\label{sdde}\\
x(t)&=&\xi(t),\quad t\in[-\de,0]\label{initial}
\end{eqnarray}
with the initial conditions $x(0)=x_0\in \mathbb{S}$.
Here
$f: \RR^n\K \RR^n\K\mathbb{S}  \to \RR^n  \quad\hbox{and}\quad g: \RR^n\K \RR^n\K\mathbb{S}  \to \RR^{n\K m}$ are measurable mapping,
 $\delta(t):[0,\infty)\rightarrow [0,\de]$ is a Borel measurable function.
We assume that the coefficients $f$ and $g$ can be decomposed as
$$
f(x,y,i)=F_1(x,y,i)+F(x,y,i)\quad \textit{and} \quad g(x,y,i)=G_1(x,y,i)+G(x,y,i).
$$

Moreover, let $C^{2}(\RR^n\times \mathbb{S};\RR_+)$ denotes the family of all nonnegative
functions $V(x,i)$ on $\RR^n\times \mathbb{S}$ which are continuously twice differentiable
in $x$. For each $V\in C^{2}(\RR^n\times \mathbb{S};\RR_+)$,
define an operate $LV$ from $C^{2}(\RR^n\times \mathbb{S};\RR_+)$ to $\RR$ by
\begin{eqnarray}\label{ito}
\nonumber
LV(x,y,i)&=&V_x(x,i)f(x,y,i)+\frac{1}{2}\trace[g^{T}(x,y,i)V_{xx}(x,i)g(x,y,i)]
+\sum_{j=1}^N\gamma_{ij} V(x,j)
\end{eqnarray}
where
\begin{eqnarray*}
 V_x(x,i)=\big(\frac{\partial V(x,i)}{\partial x_1},\ldots,\frac{\partial V(x,i)}{\partial x_n} \big),\quad
V_{xx}(x,t,i)=\big(\frac{\partial^2 V(x,i)}{\partial x_i \partial x_j}\big)_{n\times n}.
\end{eqnarray*}

The discrete Markovian chain $\{r_k^\D,k=0,1,2,\ldots\}$ can be simulated as follows:
Let $r_0^{\D}=i_0$ and generate a random number $\xi_1$ which is uniformly distributed in $[0,1]$. Define
\begin{equation*} 
 r_{1}^\Delta=
 \begin{cases}
 i_{1}& \text{if\quad $ i_{1}\in S-\{N\}$ such that$\sum\limits_{j=1}^{i_{1}-1}P_{i_{0}, j}(\Delta)\leq\xi_{1}<\sum\limits_{j=1}^{i_{1}}P_{i_{0}, j}(\Delta) $},\\
 N & \text{if\quad  $\sum\limits_{j=1}^{N-1}P_{i_{0}, j}(\Delta)\leq\xi_{1}$},
\end{cases}
\end{equation*}
where we set $\sum_{i=1}^0P_{i_0,j}(\D)=0$ as usual. Generate independently a new random number $\xi_2$ which is
again uniformly distributed in $[0,1]$ and then define
\begin{equation*} 
 r_{2}^\Delta=
 \begin{cases}
 i_{2}& \text{if\quad $ i_{2}\in S-\{N\}$ such that$\sum\limits_{j=1}^{i_{2}-1}P_{i_{0}^\Delta, j}(\Delta)\leq\xi_{2}<\sum\limits_{j=1}^{i_{2}}P_{i_{1}^\D, j}(\Delta) $},\\
 N & \text{if\quad  $\sum\limits_{j=1}^{N-1}P_{i_{1}^\Delta, j}(\Delta)\leq \xi_{2}$}.
\end{cases}
\end{equation*}

After explaining how to simulate the discrete Markovian Chain, we can now
define the partially truncated EM numerical solutions, we first choose a
strictly increasing continuous
function $\mu:\RR_+\to\RR_+$ such that $\mu(w)\to\8$ as $w\to\8$ and
\begin{equation}\label{mudef}
\sup_{|x|\ve |y|\le w} \big( |F(x,y,i)|\ve |G(x,y,i)| \big) \le \mu(w), \quad \forall w\ge 1.
\end{equation}
Denote by $\mu^{-1}$ the inverse function of $\mu$ and we see that
$\mu^{-1}$ is a strictly increasing continuous
function from $[\mu(0),\8)$ to $\RR_+$.  We also choose
a constant $\D^*\in (0,1]$ and a strictly
decreasing function $h: (0,\D^*]\to (0,\8)$ such that
\begin{equation}\label{hdef}
h(\D^*)\ge \mu(1), \ \
\lim_{\D \to 0} h(\D) = \8 \quad\hbox{and}\quad
\D^{1/4} h(\D) \le 1, \ \ \forall \D \in (0,\D^*].
\end{equation}
For a given step size $\D\in (0,\D^*]$, let us define a mapping
$\pi_\D$ from $\RR^n$ to the closed ball
$\{x\in\RR^n: |x|\le \mu^{-1}(h(\D))\}$ by
$$
\pi_\D(x) = (|x|\we \mu^{-1}(h(\D)))\, \frac{x}{|x|},
$$
where we set  $x/|x|=0$ when $x=0$.
That is, $\pi_\D$ will map $x$ to itself when $|x|\le \mu^{-1}(h(\D))$
and to $\mu^{-1}(h(\D)) x/|x|$ when $|x|> \mu^{-1}(h(\D))$.
We then define the truncated functions
\begin{equation*}
F_\D(x,y,i) = F(\pi_\D(x),\pi_\D(y),i)
\quad\hbox{and}\quad
G_\D(x,y,i) = G(\pi_\D(x),\pi_\D(y),i)
\end{equation*}
for $x,y\in\RR^n$.
It is easy to see that
\begin{equation}\label{2.6}
|F_\D(x,y,i)|\ve |G_\D(x,y,i)| \le \mu(\mu^{-1}(h(\D))) = h(\D),\quad
\forall x,y\in \RR^n.
\end{equation}
That is, both truncated functions $F_\D$ and $G_\D$ are bounded
although $F$ and $G$ may not.

Let us now form the discrete-time partially truncated EM solutions.
Define $t_k = k\D$ for $k=-M, -(M-1), \cdots, 0, 1, 2, \cdots$, where $M=\lfloor \delta(k\D)/\D\rfloor+1$.
Set $X_\D(t_k) = \xi(t_k)$ for $k=-M, -(M-1), \cdots, 0$ and then form
\begin{eqnarray} \label{TEM1}
\nonumber
& &X_{k+1} -X_k \\
\nonumber
&=&f_\D(X_k, X_{k-\lfloor\delta(k\Delta)/\Delta\rfloor},r_{k}^\Delta)\D
\nonumber
 +G_\D(X_k, X_{k-\lfloor\delta(k\Delta)/\Delta\rfloor},r_{k}^\Delta)\D B_k.\\
&=&[F_1(X_k, X_{k-\lfloor\delta(k\Delta)/\Delta\rfloor},r_{k}^\Delta)+F_\D(X_k, X_{k-\lfloor\delta(k\Delta)/\Delta\rfloor},r_{k}^\Delta)]\D \\
 \nonumber
 & &+[G_1(X_k, X_{k-\lfloor\delta(k\Delta)/\Delta\rfloor},r_{k}^\Delta)+ G_\D(X_k,X_{k-\lfloor\delta(k\Delta)/\Delta\rfloor},r_{k}^\Delta)]\D B_k.
\end{eqnarray}
For $k = 0, 1,2,\cdots$,  where $\D B_k=B(t_{k+1})-B(t_{k})$.
In our analysis, it is more convenient
to work on the continuous-time approximations.  There are two continuous-time
versions.  One is the continuous-time step process
$z_1(t), z_2(t)$ and $\bar{r}(t)$ on $t\in [-\de,\8)$ defined by
\begin{eqnarray}
\nonumber
& &z_{1}(t) = \sum_{k=0}^\8 X_k I_{[k\D, (k+1)\D)}(t),
\nonumber
\quad z_{2}(t) = \sum_{k=0}^\8 X_{k-\lfloor\delta(k\Delta)/\Delta\rfloor} I_{[k\D, (k+1)\D)}(t),\\
\nonumber
& &\bar{r}(t)=\sum_{k=0}^\8 r_{k}^\Delta I_{[k\D, (k+1)\D)}(t).
\end{eqnarray}
The other one is the continuous-time continuous process
$x_\D(t)$ on $t\in [-\de,\8)$ defined by $x_\D(t) = \xi(t)$ for $t\in [-\de,0]$
while for $t \ge 0$
\begin{eqnarray} \label{TEM3a}
\nonumber
x_\D(t)
&=&  \xi(0) + \int_0^t[F_1(z_1(s),z_2(s),\bar{r}(s))+ F_\D(z_1(s),z_2(s),\bar{r}(s))]ds \\
& &+ \int_0^t[G_1(z_1(s),z_2(s),\bar{r}(s))+ G_\D(z_1(s),z_2(s),\bar{r}(s))]dB(s).
\end{eqnarray}
We see that $x_\D(t)$ is an It\^o process on $t\ge 0$ with its It\^o differential
\begin{eqnarray} \label{TEM3}
\nonumber
dx_\D(t)
&=& [F_1(z_1(t),z_2(t),\bar{r}(t))+ F_\D(z_1(t),z_2(t),\bar{r}(t))]dt\\
& & + [G_1(z_1(t),z_2(t),\bar{r}(t)) +G_\D(z_1(t),z_2(t),\bar{r}(t))]dB(t).
\end{eqnarray}
It is useful to know that $X_k=x_\D(k\D)$ for every $k\ge -M$,
namely they coincide at $k\D$.

To analyze the partially truncated Euler-Maruyama method as well as to simulate the approximate solution,
we will need the following lemma (see \cite{Arn}). And we impose those standing hypotheses.
\begin{lemma}\label{L2.5}
Given $\D >0$ let $r_{k}^\Delta=r(k\Delta)$ for $k\geq 0$, Then $\{r_k^\Delta ,k=0,1,2,\ldots\}$
is a discrete Markov chain with the one-step transition probability matrix
\begin{equation*}
P(\Delta)=(P_{ij}(\Delta))_{N\times N}=e^{\D\Gamma}.
\end{equation*}
\end{lemma}
\begin{assp}\label{A2.1}
 There are constant $K>0$ and $\r\ge 0$ such that
 \begin{equation}\label{2.2}
|F_1(x,y,i)-F_1(\bar x,\bar y,i)|^2 \ve |G_1(x,y,i)-G_1(\bar x,\bar y,i)|^2
 \le K_1(|x-\bar x|^2 + |y-\bar y|^2)
\end{equation}
and
\begin{equation}
|F(x,y,i)-F(\bx,\by,i)|^2\ve |G(x,y,i)-G(\bx,\by,i)|^2\le K_1(1+|x|^\r+|y|^\r+|\bx|^\r+|\by|^\r)(|x-\bx|^2+|y-\by|^2)
\end{equation}
for those $x,y,\bar x, \bar y\in \RR^n$
and $ \forall i \in \mathbb{S}$.
\end{assp}

We can derive from (\ref{2.2}) that the coefficients $F_1$ and $G_1$ satisfy the linear growth condition
that there exists a constant $K_1$ such that
\begin{equation}\label{2.8}
|F_1(x,y,i)|\vee |G_1(x,y,i)|\leq K_1(1+|x|+|y|)
\end{equation}
for all $(x,y,i)\in\RR^n\K\RR^n\K \mathbb{S} $.

We can derive from Assumption \ref{A2.1} that $f$ and $g$ also satisfy
\begin{equation}\label{2.4}
|f(x,y,i)-f(\bx,\by,i)|^2\ve |g(x,y,i)-g(\bx,\by,i)|^2 \le K_1(1+|x|^\r+|y|^\r+|\bx|^\r+|\by|^\r)(|x-\bx|^2+|y-\by|^2)
\end{equation}
for those $x,y,\bar x, \bar y\in \RR^n$ and $ \forall i \in \mathbb{S}$.

\begin{assp} \label{A2.2}
There is a pair of constants $\bp\geq2$ and $K_2>0$ such that
\begin{equation}\label{KhasC}
x^T F(x,y,i) +\frac{\bp-1}{2}  |G(x,y,i)|^2 \le K_2(1+|x|^2+|y|^2)
\end{equation}
for all $(x,y,i)\in\RR^n\K\RR^n\K \mathbb{S} $.
\end{assp}

The following Lemma show that the truncated functions $F_\D$ and $G_\D$ preserve the generalized Khasminskii-type condition
for any $\D\in (0,\D^*]$ as shown Lemma 4.2 in \cite{Guo2017} and we state it here as a Lemma for the use of this paper.
\begin{lemma}
 Let Assumption \ref{A2.2} hold.
Then, for every $\D\in (0,\D^*]$ and $\bp\geq 2$, we have
\begin{equation}\label{2.7}
x^T F_\D(x,y,i) + \frac{\bp-1}{2} |G_\D(x,y,i)|^2 \le 2K_2(1+|x|^2+|y|^2)
\end{equation}
for all $(x,y,i)\in\RR^n\K\RR^n\K \mathbb{S} $.
\end{lemma}

Following a very similar approach used for (2.8) in \cite{GLMY2017}, we can show the following Lemma.
\begin{lemma}\label{L2.4}
Let Assumptions \ref{A2.1} and \ref{A2.2} hold.
Then, for every $\D\in (0,\D^*]$, we can show that for any $p\in[2,\bp)$
\begin{equation*}
x^T f(x,y,i) + \frac{p-1}{2} |g(x,y,i)|^2 \le K_3(1+|x|^2+|y|^2)
\end{equation*}
for all $(x,y,i)\in\RR^n\K\RR^n\K \mathbb{S} $.
where
$$
K_3=2K_1+K_2+\frac{K_1^2(p-1)(\bp-1)}{2(\bp-p)}.
$$
\end{lemma}

In the same way as lemma \ref{L2.4} was proved, we can show that for any $p\in[2,\bp)$,
\begin{equation}\label{2.11}
x^T f_\D(x,y,i) + \frac{p-1}{2} |g_\D(x,y,i)|^2 \le K_4(1+|x|^2+|y|^2)
\end{equation}
for all $x,y\in \RR^n$, where
$$
K_4=2K_1+2K_2+\frac{K_1^2(p-1)(\bp-1)}{2(\bp-p)}.
$$

\section{Convergence}
\noindent

We can therefore state a known result (see \cite{MY06}) as a lemma for the use of this paper.
\begin{lemma} \label{L2.3}
Let Assumptions \ref{A2.1} and \ref{A2.2} hold. Then for any given
initial data (\ref{initial}), there is a unique global solution $x(t)$ to
equation (\ref{sdde}) on $t \in [-\de,\8)$. Moreover, the solution
has the property that
\begin{equation*}\label{bd}
 \E|x(t)|^2   < \8,  \quad\forall t>0.
\end{equation*}
\end{lemma}

The following Lemma gives an upper bound, independent of $\D$, for the $p$-th moment.
\begin{lemma}\label{lemma 3.2}
Let Assumptions \ref{A2.1} and \ref{A2.2} hold. Then for any $p\in[2,\bp)$, we have
\begin{equation}\label{2.15}
\sup_{0<\D \le \D^*} \sup_{0\le t\le T} \E|x_\D(t)|^p \le C,
\end{equation}
where, and from now on, $C$ stands for generic positive real constants dependent on $T, p,\xi$ and $K_1$ etc. as well in the next sections
but independent of $\D$ and its values may change between occurrences.
\end{lemma}
\noindent
{\it Proof.} Fix $\D \in (0,\D^*]$ and the initial data $\xi$ arbitrarily.
By the general It\^o formula (\ref{ito}), we derive from (\ref{2.11}) that for $0\le t \le T$,
\begin{eqnarray}\label{3.2}
\nonumber
& &\E|x_\D(t)|^p-|\xi(0)|^p\\
\nonumber
& \leq &\E\int_0^t p|x_\D(s)|^{p-2} \Big( x_\D^T(s) f_\D(z_1(s),z_2(s),\bar{r}(s))
+\frac{p-1}{2}|g_\D(z_1(s),z_2(s),\bar{r}(s))|^2 \Big) ds \\
\nonumber
&= & \E\int_0^t p|x_\D(s)|^{p-2} \Big(z_{1}^T(s) f_\D(z_1(s),z_2(s),\bar{r}(s))
+\frac{p-1}{2} |g_\D(z_1(s),z_2(s),\bar{r}(s))|^2 \Big) ds \\
\nonumber
& &+ \E\int_0^t p|x_\D(s)|^{p-2}(x_\D(s)- z_{1}(s))^T f_\D(z_1(s),z_2(s),\bar{r}(s)) ds\\
\nonumber
&\leq & \E \int_0^t p|x_\D(s)|^{p-2}K_4 (1+z_1(s)|^2+|z_2(s)|^2) ds \\
\nonumber
& & +\E\int_0^t p|x_\D(s)|^{p-2}(x_\D(s)- z_{1}(s))^T f_\D(z_1(s),z_2(s),\bar{r}(s)) ds\\
&=& J_1+J_2,
\end{eqnarray}
where
\begin{eqnarray*}
& &J_1=\E \int_0^t p|x_\D(s)|^{p-2}K_4 (1+z_1(s)|^2+|z_2(s)|^2) ds,\\
& &J_2=\E\int_0^t p|x_\D(s)|^{p-2}(x_\D(s)- z_{1}(s))^T f_\D(z_1(s),z_2(s),\bar{r}(s)) ds.
\end{eqnarray*}
By Young inequality $a^{p-2} b \leq \frac{p-2}{p} a^p+\frac{p}{2} b^{p/2},\ \forall a,b \geq 0 $ and elementary
$(a+b+c)^p \leq 3^{p-1}(a^p+b^p +c^p)$, we then have
\begin{eqnarray}\label{3.3}
\nonumber
J_1 &=&\E \int_0^t p|x_\D(s)|^{p-2}K_4 (1+z_1(s)|^2+|z_2(s)|^2) ds\\
\nonumber
&\leq & (p-2)\E\int_0^t |x_\D(s)|^{p} ds+\frac{p}{2}\E\int_0^t (K_4 (1+z_1(s)|^2+|z_2(s)|^2))^{p/2} ds\\
\nonumber
&\leq & (p-2)\E\int_0^t |x_\D(s)|^{p} ds+\frac{p}{2}K_4^{p/2}3^{p-1} \E\int_0^t(1+z_1(s)|^p+|z_2(s)|^p) ds\\
&\leq & C\int_0^t (1+\E|x_\D(s)|^p+\E|z_1(s)|^p+\E|z_2(s)|^p) ds.
\end{eqnarray}
However, $f_\D(z_1(s),z_2(s),\bar{r}(s))=F_1(z_1(s),z_2(s),\bar{r}(s))+F_\D(z_1(s),z_2(s),\bar{r}(s)) $,
therefore,
\begin{eqnarray*}
J_2&=&\E\int_0^t p|x_\D(s)|^{p-2}(x_\D(s)- z_{1}(s))^T f_\D(z_1(s),z_2(s),\bar{r}(s)) ds\\
&=&\E\int_0^t p|x_\D(s)|^{p-2}(x_\D(s)- z_{1}(s))^T F_1(z_1(s),z_2(s),\bar{r}(s)) ds\\
& &+\E\int_0^t p|x_\D(s)|^{p-2}(x_\D(s)- z_{1}(s))^T F_\D(z_1(s),z_2(s),\bar{r}(s)) ds\\
& =&J_3+J_4,
\end{eqnarray*}
where
\begin{eqnarray*}
J_3&=&\E\int_0^t p|x_\D(s)|^{p-2}(x_\D(s)- z_{1}(s))^T F_1(z_1(s),z_2(s),\bar{r}(s)) ds,\\
J_4&=&\E\int_0^t p|x_\D(s)|^{p-2}(x_\D(s)- z_{1}(s))^T F_\D(z_1(s),z_2(s),\bar{r}(s)) ds.
\end{eqnarray*}
Similarly, by (\ref{2.8}), we also show that
\begin{eqnarray}
J_3\leq C\int_0^t (1+\E|x_\D(s)|^p+\E|z_1(s)|^p+\E|z_2(s)|^p) ds.
\end{eqnarray}
Moreover, according to (\ref{2.6}) and Young's inequality, we have
\begin{eqnarray}\label{3.5}
\nonumber
J_4& =&\E\int_0^t p|x_\D(s)|^{p-2}(x_\D(s)- z_{1}(s))^T F_\D(z_1(s),z_2(s),\bar{r}(s)) ds\\
\nonumber
&\leq & (p-2)\E\int_0^t |x_\D(s)|^p ds \\
\nonumber
& &+\frac{p}{2}\E \int_0^t|x_\D(s)-z_1(s)|^{p/2}|F_\D(z_1(s),z_2(s),\bar{r}(s))|^{p/2} ds\\
& \leq & (p-2)\E\int_0^t |x_\D(s)|^p ds +\frac{p}{2} h^{p/2}(\D) \int_0^t \E|x_\D(s)-z_1(s)|^{p/2} ds.
\end{eqnarray}
On the other hand, for any $s\in[0,T]$, there is a unique $k\geq 0$ such that $k\D \leq s <(k+1)\D$. By (\ref{2.6}), (\ref{2.8}),
element inequality and It\^o isometry, we then derive from (\ref{TEM3a}) that
\begin{eqnarray}\label{1}
\nonumber
& &\E|x_\D(t)-z_1(t)|^{p/2}=\E|x_\D(t)-x_\D(k\D)|^{p/2}\\
& \leq & 4^{p/2-1}\Big[\E|\int_{k\D}^t F_{\D}(z_1(s),z_2(s),\bar{r}(s))ds|^{p/2}
\nonumber
+\E|\int_{k\D}^t F_1(z_1(t),z_2(t),\bar{r}(t))ds|^{p/2}\\
& & +\E|\int_{k\D}^t G_{\D}(z_1(s),z_2(s),\bar{r}(s)))dB(s)|^{p/2}
\nonumber
+\E|\int_{k\D}^t G_1(z_1(s),z_2(s),\bar{r}(s)))dB(s)|^{p/2}\Big]\\
& \leq & 4^{p/2-1}\Big[\D^{p/2-1}\E \int_{k\D}^t |F_{\D}(z_1(s),z_2(s),\bar{r}(s))|^{p/2} ds
\nonumber
+\D^{p/2-1}\E\int_{k\D}^t |F_1(z_1(t),z_2(t),\bar{r}(t))|^{p/2} ds\\
& & +\D^{{p}/{4}-1}\E\int_{k\D}^t |G_{\D}(z_1(s),z_2(s),\bar{r}(s)))|^{p/2} ds
\nonumber
+\D^{p/4-1}\E\int_{k\D}^t |G_1(z_1(s),z_2(s),\bar{r}(s)))|^{p/2} ds\Big]\\
\nonumber
& \leq & 4^{p/2-1}\Big[2 \D^{p/4}h^{p/2}(\D)\\
& &+\D^{p/2-1}\E\int_{k\D}^t |F_1(z_1(t),z_2(t),\bar{r}(t))|^{p/2} ds
\nonumber
+\D^{p/4-1}\E\int_{k\D}^t |G_1(z_1(s),z_2(s),\bar{r}(s)))|^{p/2} ds\Big]\\
& \leq& C\D^{p/4} (h^{p/2}(\D)+ 1+\E|z_1(s)|^{p/2}+\E|z_2(s)|^{p/2}).
\end{eqnarray}
Substituting this into the (\ref{3.5}) and recalling (\ref{hdef}), we obtain
\begin{eqnarray}\label{3.6}
\nonumber
J_4&\leq& (p-2)\E\int_0^t |x_\D(s)|^p ds +C h^{p/2}(\D) \D^{p/4}\int_0^t(1+h^{p/2}(\D)+\E|z_1(s)|^{p/2}+\E|z_2(s)|^{p/2}) ds\\
&\leq& C(1+\int_0^t (1+\E|x_\D(s)|^p+\E|z_1(s)|^{p}+\E|z_2(s)|^{p})ds).
\end{eqnarray}
Substituting (\ref{3.3})-(\ref{3.6}) into (\ref{3.2}) yields
\begin{eqnarray*}
\E|x_\D(t)|^p \le C (1+\int_0^t  \sup_{0\le u\le s} \E|{x}_\D(u)|^p ds).
\end{eqnarray*}
As this holds for any $t\in [0,T]$, while the sum of the right-hand-side (RHS) terms
is non-decreasing in $t$, we then see
\begin{eqnarray*}
\sup_{0\le u\le t} \E|x_\D(u)|^p
\le C(1 + \int_0^t \sup_{0\le u\le s} \E|x_\D(u)|^p ds).
\end{eqnarray*}
The well-known Gronwall inequality yields that
$$
\sup_{0\le u\le T} \E|x_\D(u)|^2 \le C.
$$
As this holds for any $\D \in (0,\D^*]$, while $C$ is independent of $\D$, we obtain the
required assertion (\ref{2.15}).
\hfill $\Box$

\begin{lemma}\label{L2.6}
Let Assumption \ref{A2.1} and \ref{A2.2} hold, then for any $\D \in (0,\D^*]$, we have
\begin{equation}
\E|x_\D(t)-z_1(t)|^p \le C \D^{p/2}h^p(\D), \quad \forall t\ge 0,
\end{equation}
\end{lemma}
\noindent
{\it Proof.} By Lemma \ref{lemma 3.2}, there is a $\D\in(0,\D^*]$ such that
\begin{equation}\label{3.8}
\sup_{0<\D\leq \D^*}\sup_{0\leq t\leq T}\E|x_\D(t)|^p\leq C.
\end{equation}
Now, fix any $\D\in(0,\D^*]$, For any $t\in[0,T]$, there is a unique $k\geq 0$ such that $k\D \leq t \leq (k+1)\D$. In the same way as
(\ref{1}) was proved, we can then show
$$
\E|x_\D(t)-z_1(t)|^p \leq C\D^{p/4}(1+\E|z_1(s)|^{p/2}+\E|z_2(s)|^{p/2}+h^{p/2}(\D)).
$$
By (\ref{3.8}), we therefore have
$$
\E|x_\D(t)-z_1(t)|^p\leq C\D^{p/2}h^p(\D).
$$
The proof is complete. \hfill $\Box$

\begin{assp}\label{A3.4}
We assume that delay function $\delta(\cdot)$ is bounded and differentiable, moreover, its derivative is bounded by a constant $\overline{\delta}\in[0,1)$, that is
\begin{equation}
\frac{d\delta}{dt}\leq \overline{\delta} \quad
\textit{and}\quad
\T:=\sup_{t\geq 0}\de(t)<\infty.
\end{equation}
\end{assp}

According this Assumption \ref{A3.4}, we can obtain that there exists a constants $K_5>1$ such that
\begin{equation}\label{3}
|\delta(a)-\delta(b)|\leq K_5|a-b|.
\end{equation}
\begin{lemma}\label{4}
Let Assumption \ref{A3.4} hold, then, for any $\D \in (0,\D^*]$ and $p\geq 2$, we have
\begin{equation*}
 \E|x_{\Delta}(t-\delta(t))-z_{2}(t)|^p\leq C\Delta^{p/2} h^p(\Delta).
\end{equation*}
\end{lemma}
\noindent
{\it Proof.} In the same way to the (\ref{1}), we can obtain
\begin{eqnarray}\label{2.22}
\nonumber
& &\E|x_{\Delta}(t-\delta(t))-z_{2}(s)|^p
\nonumber
=\E|x_{\Delta}(t-\delta(t))-x_{\Delta}(k\Delta-In[\delta(k\Delta)/\Delta]\Delta)|^p\\
&\leq &C(t-\delta(t)-k\Delta-\lfloor\delta(k\Delta)/\Delta\rfloor\Delta)^{p/2}\big(1+\E|z_1(t)|^p+\E|z_2(t)|^p+h^p(\D) \big).
\end{eqnarray}
Noting
$$\delta(k\Delta)-\Delta\leq \lfloor\delta(k\Delta)/\Delta\rfloor\Delta\leq \delta(k\Delta),$$
and (\ref{3}),
we derive that
\begin{equation*}
|\delta(t)-\lfloor\delta(k\Delta)/\Delta\rfloor\Delta|\leq
 \begin{cases}
 \delta(t)-\delta(k\Delta)+\Delta \leq (K_5+1)\Delta & \text{if\quad $ \delta(t)>\delta(k\Delta)$}, \\
\delta(k\Delta)-\delta(t) \leq K_5\Delta & \text{if\quad $ \delta(t)<\delta(k\Delta)-\Delta$ },\\
\D  & \text{if\quad otherwise}.
\end{cases}
\end{equation*}
In other words, we always have
$$|\delta(t)-\lfloor\delta(k\Delta)/\Delta\rfloor\Delta|\leq (K_5+1)\Delta.$$
Therefore
$$
|t-\delta(t)-k\Delta-\lfloor\delta(k\Delta)/\Delta\rfloor\Delta|^{p/2} \leq ((K_5+2)\D)^{p/2}.
$$
Substituting this into (\ref{2.22}) gives
$$\E|x_{\Delta}(t-\delta(t))-z_{2}(s)|^p\leq C \Delta^{p/2}\big(1+\E|z_1(t)|^p+\E|z_2(t)|^p+h^p(\D) \big).$$
By Lemma \ref{lemma 3.2}, we therefore have
$$
\E|x_{\Delta}(t-\delta(t))-z_{2}(s)|^p\leq C\Delta^{p/2} h^p(\D).
$$
Then the proof is complete. \hfill $\Box$
\begin{lemma}\label{L3.7}
Let Assumption \ref{A2.1},\ref{A2.2} and \ref{A3.4} hold. For any real number $R>|x(0)|$, define the stopping time
\begin{equation*}
\T_R=\inf\{t\geq 0:|x(t)|\geq R\},
\end{equation*}
where throughout this paper we set $\inf \emptyset = \infty$ (and as usual $\emptyset$ denotes the empty set). Then
\begin{equation}\label{3.12}
\mathbb{P}(\T_{R} \leq T)\leq \frac{C}{R^p}.
\end{equation}
\end{lemma}
\noindent
{\it Proof.} By the general It\^o formula, Young inequality and Lemma \ref{L2.4}, we derive that for $0\leq t\leq T$,
\begin{eqnarray*}
& &\E|x(t\we\T_R)|^p-|\xi(0)|^p\\
&\leq & \E \int_0^{t\we\T_R}p|x(s)|^{p-2}(x^T(s) f(x(s),x(s-\delta(s)),{r}(s))\\
& &+\frac{p-1}{2}|g(x(s),x(s-\delta(s)),{r}(s))|^2) ds\\
&\leq &K_3\E \int_0^{t\we\T_R}p|x(s)|^{p-2}(1+|x(s)|^2+|x(s-\delta(s))|^2) ds\\
&\leq &C\E \int_0^{t\we\T_R}(1+|x(s)|^p+|x(t-\delta(t))|^p) ds\\
&\leq &C T+C\E\int_0^{t\we \T_R} |x(s)|^p ds+C\E\int_0^{t\we \T_R}|x(s-\delta(s))|^p ds\\
&\leq &C T+C\E\int_0^{t\we \T_R} |x(s)|^p ds
+\frac{C}{1-\bar{\delta}}\E\int_{-\delta(0)}^{(t-\delta(t))\we \T_R}|x(s)|^p ds\\
&\leq & CT+\frac{C}{1-\bar{\delta}}\E\int_{-\delta(0)}^{t\we\T_R} |x(s)|^p ds\\
&\leq &C T+\frac{C\T}{1-\bar{\delta}}|\xi|^p +\frac{C}{1-\overline{\delta}} \E\int_0^{t\we \T_R}|x(s)|^p ds\\
&\leq & C+\frac{C}{1-\bar{\delta}}\int_0^t\E|x(s\we \T_R)|^p ds.
\end{eqnarray*}
The Gronwall inequality shows
$$
\E|x(t\we \T_R)|^p \leq C.
$$
This implies, by the Chebyshev inequality,
$$
R^p\mathbb{P}(\T_R\leq T)\leq \E|x(t\we \T_R)|^p\leq C
$$
and the assertion (\ref{3.12}) follows.
\hfill $\Box$

The follows Lemma can be proved in the same way as lemma \ref{L3.7} was proved.
\begin{lemma}\label{L3.8}
let Assumption and hold. For any real number $R>|x(0)|$, define the stopping time
\begin{equation*}
\rho_{\D,R}=\inf\{t\geq 0:|x_\D(t)|\geq R\},
\end{equation*}
Then
\begin{equation}\label{2.23}
\mathbb{P}(\rho_{\D,R}\leq T)\leq \frac{C}{R^p}.
\end{equation}
\end{lemma}

In order obtain our main convergence rate theorem, we need some addition condition.
\begin{assp}\label{A3.8}
There is a pair of  constants $K_6>0$ and $v\in (0,1]$ such that the initial
data $\xi$ satisfies
$$
|\xi(a)-\xi(b)|\le K_6|a-b|^v, \quad -\T\le a<b\le 0.
$$
\end{assp}

\begin{assp}\label{A2.15}
Assume that there is a positive constant $K_7$ and $\bar{q}>2$such that
\begin{equation*}
(x-\bx)^T(F(x,y,i)-F(\bx,\by,i))+\frac{\bar {q}-1}{2}|G(x,y,i)-G(\bx,\by, i)|^2 \le K_7(|x-\bx|^2+|y-\by|^2)
\end{equation*}
for all $x,y,\bx,\by \in \RR^n$.
\end{assp}
In the same way as performed in the proof of Lemma \ref{L2.4}, and according to the Assumption \ref{A2.1}, we can shows the following Lemma.

\begin{lemma} Let Assumption \ref{A2.1} and \ref{A2.15} holds, then for any $\D\in(0,\D^*]$, we have for any $q\in (2,\bar{q})$
\begin{equation}\label{2.25}
(x-\bx)^T(f(x,y,i)-f(\bx,\by,{i}))+\frac{q-1}{2}|g(x,y,i)-g(\bx,\by,{i})|^2 \le K_7(|x-\bx|^2+|y-\by|^2)
\end{equation}
for all $x,y,\bx,\by \in \RR^n$.
\end{lemma}

The following Lemma will play a key role in the proof of the convergence rate.
\begin{lemma}\label{T3.4}
Let Assumptions \ref{A2.1}, \ref{A2.2} \ref{A3.8} and \ref{A2.15} hold and assume that ${q}> 2$ and $p\geq\r$
Let $R>|x(0)|$ be a real number and let $\D\in (0,\D^*)$ be sufficiently small
such that $\mu^{-1}(h(\D)) \ge R$.  Then   
\begin{equation}\label{3.10}
 \E|x(t\wedge\theta_{\Delta,R})-x_{\Delta}(t\wedge\theta_{\Delta,R})|^2 \le C(\D^{2 v} \ve\Delta h^2(\Delta)).
\end{equation}
\end{lemma}
\noindent
{\it Proof.} Let $\tau_R$ and $\rho_{\D,R}$ be the same as before. Let
$$
\o_{\D,R}=\T_R \we\r_{\D,R}
\quad\hbox{and} \quad
e_\D(t)= x_\D(t) - x(t).
$$
and we write $\o_{\D,R}=\o$ for simplicity.
We observe that for $0 \leq s \leq t \we \theta$,
$$
|x(s)|\vee |x_\D(s-\delta(s))| \vee |z_1(s)|\vee|z_2(s)|\leq R.
$$
Recalling the definition of the truncated
functions $F_\D$ and $G_\D$ as well as (\ref{mudef}), we hence have that
$$
F_\D(z_1(s),z_2(s),i)=F(z_1(s),z_2(s),i),
\ \
G_\D(z_1(s),z_2(s),i)=G(z_1(s),z_2(s),i)
$$
for $0\le s\le t\we\o$.
Then
$$
f_{\D}(z_1,z_2,i)=F_1(z_1,z_2,i))+F_{\D}(z_1,z_2,i)
=F_1(z_1,z_2,i)+F(z_1,z_2,i)=f(z_1,z_2,i)
$$
and
$$
g_{\D}(z_1,z_2,i)=g(z_1,z_2,i).
$$
\noindent
The It\^o formula and (\ref{2.25}) shows that
\begin{eqnarray*}\label{3.11}
& &\E|e_\D(t\wedge\theta)|^2\\
& =  &2\E\int_0 ^{t\wedge\theta} \Big((x(s)-x_\D(s))^T[ f(x(s),x(s-\delta(s)),r(s)) -f_\D(z_1(s),z_2(s),\bar{r}(s))]\\
& &+\frac{1}{2}|g(x(s),x(s-\delta(s)),r(s)) -g_\D(z_1(s),z_2(s),\bar{r}(s))|^2 \Big) ds \\
& =&  2\E\int_0 ^{t\wedge\theta} \Big((x(s)-x_\D(s))^T[ f(x(s),x(s-\delta(s)),r(s)) -f(z_1(s),z_2(s),\bar{r}(s))]\\
& &+\frac{1}{2} |g(x(s),x(s-\delta(s)),r(s)) -g(z_1(s),z_2(s),\bar{r}(s))|^2 \Big) ds \\
& =&  2\E\int_0 ^{t\wedge\theta}\Big((x(s)-z_{1}(s))^T[ f(x(s),x(s-\delta(s)),r(s)) -f(z_1(s),z_2(s),\bar{r}(s))]\\
& &+\frac{1}{2} |g(x(s),x(s-\delta(s)),r(s)) -g(z_1(s),z_2(s),\bar{r}(s))|^2 \Big) ds \\
& &+2\E\int_0 ^{t\wedge\theta} \Big (  (z_{1}(s)-x_{\Delta}(s))^T[ f(x(s),x(s-\delta(s)),r(s)) -f(z_1(s),z_2(s),\bar{r}(s))]\Big)ds \\
& =&  2\E\int_0 ^{t\wedge\theta} \Big((x(s)-z_{1}(s))^T[ f(x(s),x(s-\delta(s)),r(s)) -f(z_1(s),z_2(s),{r}(s))]\\
& &+\frac{1}{2} |g(x(s),x(s-\delta(s)),r(s)) -g(z_1(s),z_2(s),\bar{r}(s))|^2 \Big) ds \\
& &+2\E\int_0 ^{t\wedge\theta} (x(s)-z_{1}(s))^T[ f(z_1(s),z_2(s),{r}(s)) -f(z_1(s),z_2(s),\bar{r}(s))] ds\\
& &+2\E\int_0 ^{t\wedge\theta} (z_{1}(s)-x_{\Delta}(s))^T[ f(x(s),x(s-\delta(s)),r(s)) -f(z_1(s),z_2(s),\bar{r}(s))]ds \\
& =&  2\E\int_0 ^{t\wedge\theta} \Big((x(s)-z_{1}(s))^T[ f(x(s),x(s-\delta(s)),r(s)) -f(z_1(s),z_2(s),{r}(s))]\\
& &+\frac{1}{2} |g(x(s),x(s-\delta(s)),r(s)) -g(z_1(s),z_2(s),\bar{r}(s))|^2 \Big) ds \\
& &+2\E\int_0 ^{t\wedge\theta} (x(s)-z_{1}(s))^T[ f(z_1(s),z_2(s),{r}(s)) -f(z_1(s),z_2(s),\bar{r}(s))]ds\\
& &+2\E\int_0 ^{t\wedge\theta} (z_{1}(s)-x_{\Delta}(s))^T[ f(x(s),x(s-\delta(s)),r(s)) -f(z_1(s),z_2(s),{r}(s))] ds \\
& &+2\E\int_0 ^{t\wedge\theta} (z_{1}(s)-x_{\Delta}(s))^T[ f(z_1(s),z_2(s),{r}(s)) -f(z_1(s),z_2(s),\bar{r}(s))]ds .
\end{eqnarray*}
Noting
\begin{eqnarray*}
& &\frac{1}{2} |g(x(s),x(s-\delta(s)),r(s)) -g(z_1(s),z_2(s),\bar{r}(s))|^2\\
&\le& \frac{1}{2}\K\Big((1+\frac{q-2}{1})|g(x(s),x(s-\delta(s)),r(s))-g(z_1(s),z_2(s),{r}(s))|^2\\
& & +(1+\frac{1}{q-2})|g(z_1(s),z_2(s),{r}(s))-g(z_1(s),z_2(s),\bar{r}(s))|^2\Big)\\
&=&\frac{q-1}{2}|g(x(s),x(s-\delta(s)),r(s))-g(z_1(s),z_2(s),{r}(s))|^2\\
& &+\frac{q-1}{2(q-2)}|g(z_1(s),z_2(s),{r}(s))-g(z_1(s),z_2(s),\bar{r}(s))|^2.
\end{eqnarray*}
Therefore
\begin{eqnarray}\label{H}
\nonumber
& &\E|e_\D(t\wedge\theta)|^2\\
\nonumber
& \leq& 2\E\int_0 ^{t\wedge\theta} \Big((x(s)-z_{1}(s))^T[ f(x(s),x(s-\delta(s)),r(s)) -f(z_1(s),z_2(s),{r}(s))]\\
\nonumber
& &+\frac{q-1}{2} |g(x(s),x(s-\delta(s)),r(s)) -g(z_1(s),z_2(s),{r}(s))|^2 \Big) ds \\
\nonumber
& &+2\E\int_0 ^{t\wedge\theta} (x(s)-z_{1}(s))^T[ f(z_1(s),z_2(s),{r}(s)) -f(z_1(s),z_2(s),\bar{r}(s))] ds\\
\nonumber
& &+2\E\int_0 ^{t\wedge\theta} (z_{1}(s)-x_{\Delta}(s))^T[ f(x(s),x(s-\delta(s)),r(s)) -f(z_1(s),z_2(s),{r}(s))] ds \\
\nonumber
& &+2\E\int_0 ^{t\wedge\theta} (z_{1}(s)-x_{\Delta}(s))^T[ f(z_1(s),z_2(s),{r}(s)) -f(z_1(s),z_2(s),\bar{r}(s))]ds \\
\nonumber
& &+\E\int_0^{t\wedge\theta}\frac{q-1}{q-2}|g(z_1(s),z_2(s),r(s)) -g(z_1(s),z_2(s),\bar{r}(s))|^2 \Big) ds \\
&=& H_1+H_2+H_3+H_4+H_5,
\end{eqnarray}
where
\begin{eqnarray*}
H_1&=& 2\E\int_0 ^{t\wedge\theta} \Big((x(s)-z_{1}(s))^T[ f(x(s),x(s-\delta(s)),r(s)) -f(z_1(s),z_2(s),{r}(s))]\\
& &+\frac{q-1}{2} |g(x(s),x(s-\delta(s)),r(s)) -g(z_1(s),z_2(s),{r}(s))|^2 \Big) ds, \\
H_2&=&2\E\int_0 ^{t\wedge\theta} (x(s)-z_{1}(s))^T[ f(z_1(s),z_2(s),{r}(s)) -f(z_1(s),z_2(s),\bar{r}(s))] ds, \\
H_3&=&2\E\int_0 ^{t\wedge\theta} (z_{1}(s)-x_{\Delta}(s))^T[ f(x(s),x(s-\delta(s)),r(s)) -f(z_1(s),z_2(s),{r}(s))] ds, \\
H_4&=&2\E\int_0 ^{t\wedge\theta} (z_{1}(s)-x_{\Delta}(s))^T[ f(z_1(s),z_2(s),{r}(s)) -f(z_1(s),z_2(s),\bar{r}(s))]ds
\end{eqnarray*}
and
\begin{eqnarray*}
H_5&=&\E\int_0^{t\wedge\theta}\frac{q-1}{q-2}|g(z_1(s),z_2(s),r(s)) -g(z_1(s),z_2(s),\bar{r}(s))|^2  ds\ \ \ \ \
\end{eqnarray*}
According to the Young's inequality and the (\ref{2.25}), we have
\begin{eqnarray*}
H_1&=& 2\E\int_0 ^{t\wedge\theta} \Big((x(s)-z_{1}(s))^T[ f(x(s),x(s-\delta(s)),r(s)) -f(z_1(s),z_2(s),{r}(s))]\\
& &+\frac{q-1}{2} |g(x(s),x(s-\delta(s)),r(s)) -g(z_1(s),z_2(s),{r}(s))|^2 \Big) ds \\
&\leq& 2K_{7}\E\int_0 ^{t\wedge\theta}(|x(s)-z_{1}(s)|^2+|x(s-\delta(s))-z_{2}(s)|^2)ds\\
&\leq &2K_{7}\E\int_0 ^{t\wedge\theta}\big (|x(s)-x_{\Delta}(s)|^2+|x_{\Delta}(s)-z_{1}(s)|^2\\
& &+|x(s-\delta(s))-x_{\Delta}(s-\delta(s))|^2+|x_{\Delta}(s-\delta(s))-z_{2}(s)|^2)ds\\
&\leq & 2K_{7}\E\int_0 ^{t\wedge\theta} |x(s)-x_{\Delta}(s)|^2 ds
+2K_{7}\int_0 ^{T}\E |x_{\Delta}(s)-z_{1}(s)|^2 ds\\
& &+2K_{7}\int_0 ^{T}\E |x_{\Delta}(s-\delta(s))-z_{2}(s)|^2 ds
+2K_{7}\E\int_{-\delta} ^{0} |\xi(\lfloor s/\D\rfloor\D)-\xi(s)|^2 ds .
\end{eqnarray*}
By Lemma \ref{L2.6}, \ref{4} and Assumption \ref{A5.2}, we have
\begin{eqnarray}\label{H_1}
H_1 \leq C \big(\int_0^t \E|x(s\we \o)-x_{\Delta}(s\we \o)|^2 ds + \D h^2(\D)+\D^{2 v}\big).
\end{eqnarray}
By the Young's inequality and elementary inequality, we can obtain
\begin{eqnarray}\label{3.28}
\nonumber
H_2&=&2\E\int_0 ^{t\wedge\theta} (x(s)-z_{1}(s))^T[ f(z_1(s),z_2(s),{r}(s)) -f(z_1(s),z_2(s),\bar{r}(s))] ds\\
\nonumber
&\leq & \E\int_0 ^{t\wedge\theta} |x(s)-z_1(s)|^2 ds+\E\int_0 ^{t\wedge\theta} | f(z_1(s),z_2(s),{r}(s)) -f(z_1(s),z_2(s),\bar{r}(s))|^2 ds\\
\nonumber
&= & \E\int_0 ^{t\wedge\theta} |x(s)-x_\D(s)+x_\D(s)-z_1(s)|^2 ds\\
\nonumber
& &+\E\int_0 ^{t\wedge\theta} | f(z_1(s),z_2(s),{r}(s)) -f(z_1(s),z_2(s),\bar{r}(s))|^2 ds\\
&\leq & \E\int_0 ^{t\wedge\theta}|x_\D(s)-x(s)|^2 ds
\nonumber
+\E\int_0 ^{t\wedge\theta} |x_\D(s)-z_1(s)|^2 ds\\
\nonumber
& &+\E\int_0 ^{t\wedge\theta} | f(z_1(s),z_2(s),{r}(s)) -f(z_1(s),z_2(s),\bar{r}(s))|^2 ds\\
&\leq & \E\int_0 ^{t\wedge\theta}|x_\D(s)-x(s)|^q ds
\nonumber
+\E\int_0 ^{t\wedge\theta} |x_\D(s)-z_1(s)|^q ds\\
& &+\E\int_0 ^{T} |f(z_{1}(s),z_{2}(s)),r(s))-f(z_{1}(s),z_{2}(s)),\bar{r}(s))|^2 ds.
\end{eqnarray}
Let $j$ be the integer part of $T/\D$. Then
\begin{eqnarray}\label{0}
\nonumber
& &\E\int_0 ^{T} |f(z_{1}(s),z_{2}(s)),r(s))-f(z_{1}(s),z_{2}(s)),\bar{r}(s))|^2 ds\\
\nonumber
&=& \sum_{k=0}^j \E\int_{k\Delta} ^{(k+1)\Delta}|f(z_1(s),z_2(s),r(s))-f(z_1(s),z_2(s),r(k\Delta))|^2 ds\\
\nonumber
&\leq &2\sum_{k=0}^j \E\int_{k\Delta} ^{(k+1)\Delta}[|f(z_1(s),z_2(s),r(s))|^2+|f(z_1(s),z_2(s),r(k\Delta))|^2] I_{\{r(s)\neq r(k\Delta)\}}ds\\
&\leq & 2\sum_{k=0}^j\int_{k\Delta} ^{(k+1)\Delta} \E[\E[(1+|z_1(s)|^2+|z_2(s)|^2 +h^2(\D))I_{\{r(s)\neq r(k\Delta)\}}|r(k\D)]]ds.
\end{eqnarray}
where in the last step, we use the fact that $z_1(s)$ and $z_2(s)$ are conditionally independent of
$I_{\{r(s)\neq r(k\Delta)\}}$ given the $\sigma-$algebra generated by $r(k\Delta)$.
But, by the Markov property
\begin{eqnarray*}
& &\E[I_{r(s)\neq r(k\Delta)}|r(k\Delta)]\\
&=&\sum_{i\in S}I_{\{r(k\Delta)=i\}}P(r(s)\neq i|r(k\Delta)=i)\\
&=&\sum_{i\in S}I_{\{r(k\Delta)=i\}}\sum_{j\neq i}(\gamma_{ij}(s-t_{k})+o(s-t_{k}))\\
&\leq &(\max_{i\leq i \leq N}(-\gamma_{ii})\Delta+o(\Delta))\sum_{i\in \mathbb{S}}I_{\{r(k\Delta)=i\}}\\
&\leq &C\Delta+o(\Delta).
\end{eqnarray*}
So, by lemma \ref{lemma 3.2},
\begin{eqnarray*}
& &\E\int_{k\Delta} ^{(k+1)\Delta}|f(z_1(s),z_2(s),r(s))-f(z_1(s),z_2(s),r(k\Delta))|^2ds\\
&\leq&(C \Delta+o(\Delta))\int_{k\Delta} ^{(k+1)\Delta} [1+\E|z_1(s)|^2+\E|z_2(s)|^2+h^2(\D)] ds\\
&\leq& h^2(\D)\D(C{\Delta}+o(\Delta)).
\end{eqnarray*}
Substituting this into (\ref{0}) gives
\begin{eqnarray*}
& &\E\int_0 ^{T} |f(z_{1}(s),z_{2}(s)),r(s))-f(z_{1}(s),z_{2}(s)),\bar{r}(s))|^2 ds\\
&\leq& h^2(\D)(C{\Delta}+o(\Delta)).
\end{eqnarray*}
This implies that
\begin{eqnarray*}
& &\E\int_0 ^{t\wedge\theta}|f(z_{1}(s),z_{2}(s),r(s))-f(z_{1}(s),z_{2}(s),r(k\Delta))|^2 ds\\
&\leq & h^2(\D)(C{\Delta}+o(\Delta)).
\end{eqnarray*}
Substituting this into (\ref{3.28}) and lemma \ref{L2.6}, we have
\begin{eqnarray}\label{H_2}
H_2 &\leq & C\big(\int_0 ^{t}\E|x(s\wedge\theta)-x_\D(s\wedge\theta)|^2 ds+\Delta h^2(\Delta)+o(\Delta)).
\end{eqnarray}
Moreover, using the (\ref{2.4}) and Lemma \ref{L2.6}, we have
\begin{eqnarray*}
H_3&=& 2 \E\int_0^{t\wedge\theta} (z_{1}(s)-x_{\Delta}(s))^T[f(x(s),x(s-\delta(s)),r(s))-f(z_1(s),z_2(s),{r}(s))] ds \\
&\leq&\E\int_0^{t\wedge\theta} |z_{1}(s)-x_{\Delta}(s)|^{2}
+\E\int_0^{t\wedge\theta}|f(x(s),x(s-\delta(s)),r(s))-f(z_1(s),z_2(s),{r}(s))|^{2} ds\\
&\leq&\E\int_0^{t\wedge\theta} |z_{1}(s)-x_{\Delta}(s)|^{2}
+ K_1\E\int_0^{t\wedge\theta}[1+|x(s)|^{\r}+|x(s-\delta(s))|^{\r}+|z_1(s)|^{\r}+|z_2(s)|^{\r}]\times\\
& &\quad (|x(s)-z_1(s)|^{2}+|x(s-\delta(s))-z_2(s)|^{2}) ds\\
&\leq &C\Big(\E\int_0^{t}|x(s\wedge\theta)-x_\D(s\wedge\theta)|^{2} ds+ \D h^2(\D)+\D^{2 v}\Big).
\end{eqnarray*}
Therefore
\begin{eqnarray}\label{H_3}
H_3 \leq C(\int_0^{t}\E|x_\D(s\wedge\theta)-x(s\wedge\theta)|^{2} ds +\D h^2(\D)+\D^{2 v}).
\end{eqnarray}
Similarly to $H_2$, we can show
\begin{eqnarray}\label{H_4}
\nonumber
H_4&=& 2\E\int_0 ^{t\wedge\theta} (z_{1}(s)-x_{\Delta}(s))^T[ f(z_1(s),z_2(s),{r}(s)) -f(z_1(s),z_2(s),\bar{r}(s))] ds \\
&\leq & C\Big(\int_0 ^{t}\E|x(s\wedge\theta)-x_\D(s\wedge\theta)|^2 ds+\Delta h^2(\Delta)+o(\Delta)\Big).
\end{eqnarray}
And
\begin{eqnarray}\label{H_5}
\nonumber
H_5 &=&\E\int_0^{t\wedge\theta}\frac{q-1}{q-2}|g(z_1(s),z_2(s),r(s)) -g(z_1(s),z_2(s),\bar{r}(s))|^2\\
&\leq & C\big(\int_0 ^{t}\E|x(s\wedge\theta)-x_\D(s\wedge\theta)|^2 ds+\Delta h^2(\Delta)+o(\Delta)).
\end{eqnarray}
Substituting (\ref{H_1}),(\ref{H_2}),(\ref{H_3}),(\ref{H_4}) and (\ref{H_5}) into (\ref{H}), we obtain that
\begin{eqnarray*}
\E|x(t\wedge\theta)-x_{\Delta}(t\wedge\theta)|^2
\leq C\Big(\int_{0}^{t}\E|x(s\wedge\theta)-x_{\Delta}(s\wedge\theta)|^2 ds +(\Delta h^2(\Delta)\ve \D^{2 v})\Big).
\end{eqnarray*}
By the well-known Gronwall inequality yields that
$$
\E|x(t\wedge\theta)-x_{\Delta}(t\wedge\theta)|^2 \leq C(\Delta h^2(\Delta)\ve \D^{2 v}).
$$
Then the proof is complete.
\hfill $\Box$

Let us now state our first result on the convergence rate.
\begin{theorem}
Let Assumptions \ref{A2.1}, \ref{A2.2}, \ref{A5.2} and \ref{A2.15} hold and assume that $p\in(2,\bp)$ and $p>\r$,
If
\begin{equation}\label{3.21}
h(\D)\geq \mu\Big((\D^{2 v} \ve \D h^2(\D))^{-1/{(p-2)}}\Big),
\end{equation}
then there is a $\D\in (0,\D^*)$ such that
\begin{equation}\label{5}
\E|x_\D(T)-x(T)|^2 \leq C(\D^{2 v}\ve \D h^{2}(\D)).
\end{equation}
\end{theorem}
\noindent
{\it Proof.} Let $\e >0$ be arbitrary. Let $\T_R, \rho_{\D,R}, \o_{\D,R}$ and $e_\D(T)$ be same as before.
For a sufficiently large $R>|x(0)|$, we have that
\begin{eqnarray}\label{3.22}
\nonumber
\E|e_\D(T)|^2 &=& \E(|e_\D(T)|^2 I_{\o_{\D,R}>T})+\E(|e_\D(T)|^2 I_{\o_{\D,R}\leq T})\\
&\leq & \E(|e_\D(T)|^2 I_{\o_{\D,R}>T})+\frac{2 \epsilon}{p}\E|e_\D(T)|^2+\frac{p-2}{p\epsilon^{2/(p-2)}}\mathbb{P}(\o_{\D,R\leq T}).
\end{eqnarray}
Applying Lemma \ref{L2.3} and \ref{lemma 3.2}, we can see that
\begin{equation}\label{3.23}
\E|e_{\D}(T)|^2 \leq 2 \E|x(T)|^2+2 \E|x_\D(T)|^2 \leq C.
\end{equation}
Using the Lemma \ref{L3.7} and \ref{L3.8}, we obtain that
\begin{equation}\label{3.24}
\mathbb{P}(\o_{\D,R}\leq T)\leq \mathbb{P}(\T_R \leq {T})+\mathbb{P}(\rho_{\D,R}\leq T) \leq \frac{C}{R^2}.
\end{equation}
Substituting (\ref{3.23}) and (\ref{3.24}) into (\ref{3.22}) and choosing $\e=\D^{2 v} \ve \D h^{2}(\D)$
and $R=(\D^{2 v} \ve \D h^{2}(\D))^{-1/{(p-2)}}$, we have that
\begin{equation}\label{3.26}
\E(|e_\D(T)|^2 I_{\o_{\D,R\leq T}}) \leq C(\D^{2 v} \ve \D h^{2}(\D)).
\end{equation}
By the Lemma \ref{T3.4}, we can show that
\begin{equation}\label{3.27}
\E(|e_{\D}(T\we \o_{\D,R})|^2) \leq C(\D^{2 v} \ve \D h^{2}(\D)).
\end{equation}
By the (\ref{3.21}), we can see that
$$
\mu^{-1}(h(\D)) \geq (\D^{2 v} \ve \D h^{2}(\D))^{-1/{ (p-2)}}=R.
$$
Therefore, substituting (\ref{3.26}) and (\ref{3.27}) into (\ref{3.22}) yields (\ref{5}).
The proof is therefore complete.
\hfill $\Box$

\begin{expl}
Consider a nonlinear scalar hybrid SDDE
\begin{equation}\label{example1}
dx(t)=f(x(t),x(t-\de(t)),r(t))dt+ g(x(t),x(t-\de(t)),r(t))dB(t)
\end{equation}
Here, $B(t)$ is a scalar Brownian, delay function $\de(t)=0.1\cos(t)$, and $r(t)$ is a Markovian chain on the state space $\mathbb{S}=\{1,2\}$ and they are independent.
Let the generator of the Markovian chain that
\begin{equation*}\Gamma=       
\left(                 
  \begin{array}{cc}   
    -2 & 2 \\  
    1 & -1   
  \end{array}
\right).                 
\end{equation*}
Moreover, for $\forall (x,y,i)\in \RR \times \RR \times \mathbb{S}$,
\begin{equation*}
f(x,y,i)=
\begin{cases}
-6x-x^5+y &\quad \text{if} \quad i=1\\
-6x-x^5+\frac{y}{1+y^2}& \quad \text{if} \quad i=2
\end{cases}
\quad and \quad
g(x,y,i)=
\begin{cases}
x^2 &\quad \text{if} \quad i=1\\
\sin x\sin y& \quad \text{if} \quad i=2
\end{cases}
\end{equation*}
\end{expl}
\noindent
\textit{Step 1.} Check the assumptions

It can be seen that
\begin{equation*}
F_1(x,y,i)=
\begin{cases}
-6x+y &\quad \text{if} \quad i=1\\
-6x& \quad \text{if} \quad i=2
\end{cases}
\quad and \quad
G_1(x,y,i)=
\begin{cases}
0 &\quad \text{if} \quad i=1\\
0 & \quad \text{if} \quad i=2
\end{cases}
\end{equation*}
\begin{equation*}
F(x,y,i)=
\begin{cases}
-x^5 &\quad \text{if} \quad i=1\\
-x^5+\frac{y}{1+y^2}& \quad \text{if} \quad i=2
\end{cases}
\quad and \quad
G(x,y,i)=
\begin{cases}
x^2 &\quad \text{if} \quad i=1\\
\sin x\sin^2y& \quad \text{if} \quad i=2
\end{cases}
\end{equation*}
then Assumption \ref{A2.1} holds and the delay function $\delta(t)=0.1\cos t$ fulfilled Assumption \ref{A3.4}, clearly. For Assumption \ref{A2.15}, it is straightforward to see that
\begin{eqnarray*}
& &(x-\bx)^T (F(x,y,1)-F(\bx,\by,1))+\frac{q-1}{2}|G(x,y,1)-G(\bx,\by,1)|^2\\
&=& (x-\bx)(-x^5+\bx^5)+\frac{q-1}{2}|x^2-\bx^2|^2\\
&=& (x-\bx)[-(x-\bx)(x^4+x^3\bx+x^2\bx^2+x\bx^3+\bx^4)]+\frac{q-1}{2}|(x-\bx)(x+\bx)|^2\\
&=& (x-\bx)^2[-(x^4+x^3\bx+x^2\bx^2+x\bx^3+\bx^4)+\frac{q-1}{2}(x+\bx)^2].
\end{eqnarray*}
However
\begin{eqnarray*}
-(x^3\bx+x\bx^3)=-x\bx(x^2+\bx^2)\leq 0.5(x^2+\bx^2)^2=0.5(x^4+\bx^4)+x^2\bx^2.
\end{eqnarray*}
Hence
\begin{eqnarray*}
& &(x-\bx)^T (F(x,y,1)-F(\bx,\by,1))+\frac{q-1}{2}|G(x,y,1)-G(\bx,\by,1)|^2\\
&\leq& (x-\bx)^2[-0.5(x^4+\bx^4)+\frac{q-1}{2}(x^2+\bx^2)]\\
&\leq &[1+\frac{(q-1)^2}{4}](x-\bx)^2.
\end{eqnarray*}
Moreover,
\begin{eqnarray*}
& &(x-\bx)^T (F(x,y,2)-F(\bx,\by,2))+\frac{q-1}{2}|G(x,y,2)-G(\bx,\by,2)|^2\\
&=& (x-\bx)[-x^5+\frac{y}{1+y^2}-(-\bx^5+\frac{\by}{1+\by^2})]+\frac{q-1}{2}|\sin x\sin y-\sin\bx \sin\by|^2.
\end{eqnarray*}
Using the mean value theorem, we set $A(x)=\frac{x}{1+x^2}$, then there is existing a $ x^* \in (x,\bx)$ such that
\begin{eqnarray*}
|A(x)-A(\bx)|=|A'(\xi)(x-\bx)|\leq |A'(x^*)||x-\bx|,
\end{eqnarray*}
where $A'(x)$ is the derivative and
$A'(x)=\frac{1-x^2}{(1+x^2)^2}\leq 1$
therefore
\begin{eqnarray*}
|\frac{y}{1+y^2}-\frac{\by}{1+\by^2}|\leq |y-\by|.
\end{eqnarray*}
Meanwhile,
\begin{eqnarray*}
& &|\sin x \sin y-\sin \bx \sin \by|^2=|\sin x \sin y-\sin x \sin \by +\sin x \sin \by-\sin \bx \sin \by|^2\\
&\leq & | \sin x \sin y-\sin x \sin \by|^2+|\sin x \sin \by-\sin \bx \sin \by|^2\\
&\leq & |y-\by|^2+|x-\bx|^2.
\end{eqnarray*}
Therefore
\begin{eqnarray*}
& &(x-\bx)^T (F(x,y,2)-F(\bx,\by,2))+\frac{q-1}{2}|G(x,y,2)-G(\bx,\by,2)|^2\\
& \leq & (x-\bx)(-0.5(x^4+\bx^4)+|y-\by|)+\frac{q-1}{2} (|x-\bx|^2+|y-\by|^2)\\
& \leq & (x-\bx)(|y-\by|)+\frac{q-1}{2}(|x-\bx|^2+|y-\by|^2)\\
&\leq & \frac{q}{2}|x-\bx|^2 +\frac{q}{2}|y-\by|^2.
\end{eqnarray*}
So, for any $i\in\mathbb{S}$, we have
\begin{equation*}
(x-\bx)^T (F(x,y,i)-F(\bx,\by,i))+\frac{q-1}{2}|G(x,y,i)-G(\bx,\by,i)|^2 \leq (\frac{q^2+3}{4})(|x-\bx|^2+|y-\by|^2).
\end{equation*}
In other words, Assumption \ref{A2.15} is also fulfilled for any $q$. Moreover,
\begin{equation*}
\begin{aligned}
 xF(x,y,1)+\frac{\bp-1}{2}|G(x,y,1)|^2
=&-x^6+\frac{\bp-1}{2} |x^2|^2\\
=&-x^2(x^2-\frac{\bp-1}{4})^2+\frac{(\bp-1)^2}{16}x^2
\leq \frac{(\bp-1)^2}{16} x^2
\end{aligned}
\end{equation*}
and
\begin{equation*}
\begin{aligned}
xF(x,y,2)+\frac{\bp-1}{2}|G(x,y,2)|^2
&=-x^6+\frac{xy}{1+y^2}+\frac{\bp-1}{2}|\sin x\sin y|^2\\
&\le -x^6+\frac{xy}{1+y^2}+\frac{\bp-1}{2}\\
&\le -x^6+ \frac{x^2}{2(1+y^2)}+ \frac{y^2}{2(1+y^2)}+\frac{\bp-1}{2}\\
& \le \frac{|x|^2}{2}+\frac{|y|^2}{2}+\frac{\bp-1}{2}\\
& \le \frac{\bp-1}{2}(1+|x|^2+|y|^2).
\end{aligned}
\end{equation*}
Therefore  for any $i\in \mathbb{S}$, we have
\begin{equation}
xF(x,y,i)+\frac{\bp-1}{2}|G(x,y,i)|^2 \leq \frac{(\bp-1)^2}{2}(1+|x|^2+|y|^2)
\end{equation}
that is, Assumption \ref{A2.2} is satisfied for any $\bp$.

\noindent
\textit{Step 2.} We need choose $\mu(\cdot)$ and $h(\cdot)$.

According to (\ref{mudef}),
$|F(x,y,i)|=|-x^5+\frac{y}{1+y^2}|\leq |x|^5$ and
$|G(x,y,i)|=|\sin x\sin y|\leq 1$, then
 we can set $\mu(w)=w^5$ such that
$$
\sup_{|x|\ve|y|\leq w}(|F(x,y,i)|\ve|G(x,y,i)|)\leq \sup_{|x|\ve|y|\leq w}|x|^5<w^5,\quad w>1.
$$
If we set $h(\D)=\D^{-1/10}$, then all the conditions in (\ref{hdef}) hold for all $\D^*\in(0,1]$, and obviously we have $\mu^{-1}(h(\D))=\D^{-1/50}$.

\noindent
\textit{Step 3.}
Applying (\ref{TEM1}), we can obtain the numerical solution. Since it is hardly to find a true solution of (\ref{example1}).
We use the numerical solution produced by partially truncated EM method with step size  $2\times10^{-4}$, $2^2\times10^{-4}$, $2^3\times10^{-4}$ and $2^4\times10^{-4}$
and the step size $10^{-4}$ as the `true solution' at the terminal time $T=1$, the square roots of the mean square errors are plotted in Figure 1 by simulating 500 paths.
We can see that the convergence rate is approximately $1/2$.
\begin{figure}
	\centering
	\includegraphics[height=6cm, width=8cm]{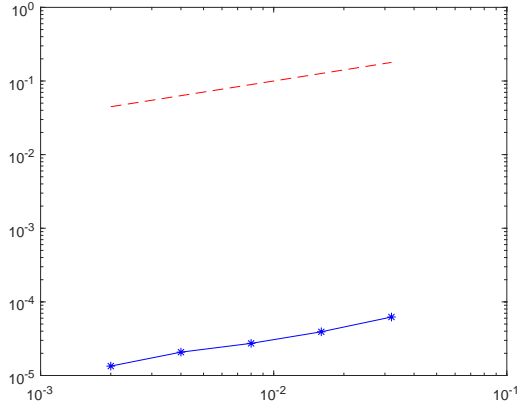}
	\caption{convergence rate}
\end{figure}
\section{Stability}
\begin{assp}\label{A5.2}
Assume that there are constants $\o\in[0,\8]$ and there is a pair of constant $\l_1>2\l_3\geq 2(1-\bar{\de})\l_4\geq 0,\l_2\geq 0$ and $\bar{\de}\in(0,1)$ such that
\begin{equation}\label{4.1}
2x^TF_1(x,y,i)+(1+\o)|G_1(x,y,i)|^2\le -\l_1|x|^2+\l_2(1-\bar{\de})|y|^2
\end{equation}
and
\begin{equation}\label{4.2}
2x^TF(x,y,i)+(1+\o^{-1})|G(x,y,i)|^2\le \l_3|x|^2+\l_4(1-\bar{\de})|y|^2
\end{equation}
for all $x,y \in \RR^d$, where throughout the remaining part of this paper.
We choose $\o=0$ and set $\o^{-1}|G(x,y)|^2=0$ when there is no $G(x,y)$ term in $g(x,y)$,
where choose $\o=\8$ and set $\o|G(x,y)|^2=0$ when there is no $G_1(x,y)$ term in $g(x,y)$.
\end{assp}

This can implying that
\begin{equation}
2x^T f(x,y,i)+|g(x,y,i)|^2 \leq -(\l_1-\l_3)|x|^2 +(1-\bar{\de})(\l_2+\l_4)|y|^2
\end{equation}
\begin{lemma}
Let Assumption \ref{A5.2} hold, then for all $\Delta\in(0,\Delta^*]$
\begin{equation}\label{fg}
2x^T f_\D(x,y,i)+|g_\D(x,y,i)|^2 \leq -(\l_1-2\l_3)|x|^2 +(1-\bar{\de})(\l_2+\l_4)|y|^2.
\end{equation}
\end{lemma}
The proof is similar to the (4.8) on \cite{GLMY2017}, so we omit it here.

From \cite{21} we known that the SDDEs (\ref{SDDE}) is almost sure exponential stable. To be precise,
we cite it as follows.
\begin{lemma}\label{L1.3}
Let Assumption \ref{A3.4} and \ref{A5.2} holds.
Then for any given initial data (\ref{initial}), the solution has the properties that
\begin{equation}
\limsup_{t\rightarrow\infty}\frac{1}{t}\log\E|x(t)|\leq-\frac{\eta }{2}
\end{equation}
and
\begin{equation}
\limsup_{t\rightarrow\infty}\frac{1}{t}\log|x(t)|\leq-\frac{\eta }{2} \quad a.s.
\end{equation}
where $\eta$ is the unique root to the equation
\begin{equation}
\l_1-2\l_3=\eta +(\l_2+\l_4)e^{\eta \T}.
\end{equation}
\end{lemma}

The following theorem shows that the partially truncated EM method can preserve this almost sure exponential stability.
\begin{theorem}\label{T4.4}
Let Assumption \ref{A3.4} and \ref{A5.2} holds.
Let $\g^*$ be the unique positive root of the equation
\begin{equation}\label{14}
[(1-\bar{\de})(\l_2+\l_3)+\e] e^{\g \T}=\l_1-2\l_3-\e-\g
\end{equation}
with a positive number $\e$ satisfied
\begin{equation}\label{12}
\e <[(1-\bar{\de})(\l_2+\l_4)-\l_1+2\l_3]/2.
\end{equation}
Then the solution generated by the partially truncated EM method is almost sure exponential stable, i.e.,
\begin{equation}
\limsup_{k\rightarrow \infty} \frac{1}{k\D} \log |X_k| \le \frac{\g}{2}+\e \ \ \ a.s.
\end{equation}
\end{theorem}
\noindent
{\it Proof .} For any positive constant $C>1$, we have
\begin{eqnarray*}
C^{(k+1)\D}|X_{k+1}|^2-C^{k\D}|X_k|^2=C^{(k+1)\D}(|X_{k+1}|^2-|X_k|^2)+(C^{(k+1)\D}-C^{k\D})|X_k|^2.
\end{eqnarray*}
For simple, we let
$$f_\D=f_\D(X_k,X_{k-\lfloor\delta(k\D)/\D\rfloor},r_{k}^\Delta)\quad\text{and}\quad
g_\D=g_\D(X_k,X_{k-\lfloor\delta(k\D)/\D\rfloor},r_{k}^\Delta).
$$
We can easily obtain from (\ref{TEM1}) that
\begin{eqnarray*}
|X_{k+1}|^2
=|X_k|^2+2 X_k^T f_\D\D+|g_\D|^2\D+|f_\D|^2\D^2+m_k.
\end{eqnarray*}
where
\begin{eqnarray*}
m_k=2 X_k^T g_{\Delta}\Delta B_k+2 \Delta f_{\Delta} g_{\Delta}\Delta B_k+|g_{\Delta}|^2((\Delta B_k)^2 - \Delta)
\end{eqnarray*}
which is a martingale (see \cite{wu2010}).
By condition (\ref{fg}) that
\begin{eqnarray}\label{15}
\nonumber
& &C^{(k+1)\D}|X_{k+1}|^2-C^{k\D}|X_k|^2\\
&=& C^{(k+1)\D}(2X_k^T f_\D \D+|g_\D|^2\D +|f_\D|^2\D^2+m_k)
\nonumber
+(C^{(k+1)\D}-C^{k\D})|X_k|^2\\
\nonumber
& \leq &-(\l_1-2\l_3)C^{(k+1)\D}\D|X_k|^2+C^{(k+1)\D}(1-\bar{\de})(\l_2+\l_4)\D|X_{k-\lfloor\delta(k\D)/\D\rfloor}|^2\\
&&+C^{(k+1)\D}|f_\D|^2\D^2+C^{(k+1)\D}m_k
\nonumber
+(C^{(k+1)\D}-C^{k\D})|X_k|^2\\
\nonumber
&\leq &-(\l_1-2\l_3)C^{(k+1)\D}\D|X_k|^2+C^{(k+1)\D}(1-\bar{\de})(\l_2+\l_4)\D|X_{k-\lfloor\delta(k\D)/\D\rfloor}|^2\\
& &+C^{(k+1)\D}\e\D(|X_k|^2+|X_{k-\lfloor\delta(k\D)/\D\rfloor}|^2)+C^{(k+1)\D} m_k
\nonumber
+(C^{(k+1)\D}-C^{k\D})|X_k|^2\\
\nonumber
&\leq & C^{(k+1)\D}(-(\l_1-2\l_3-\e)\D)|X_k|^2+C^{(k+1)\D}((1-\bar{\de})(\l_2+\l_4)\D+\e\D)|X_{k-\lfloor\delta(k\D)/\D\rfloor}|^2\\
& &+(C^{(k+1)\D}-C^{k\D})|X_k|^2+C^{(k+1)\D} m_k.
\end{eqnarray}
Applying induction to (\ref{15}) gives
\begin{eqnarray}\label{11}
\nonumber
& &C^{k\D}|X_k|^2-|X_0|^2\\
&\leq &
\nonumber
[-(\l_1-2\l_3-\e)\D]\sum_{i=0}^{k-1}C^{(i+1)\D}|X_i|^2\\
\nonumber
& &+[(1-\bar{\de})(\l_2+\l_4)\D+\e\D]\sum_{i=0}^{k-1}C^{(i+1)\D}|X_{i-\lfloor\delta(i\D)/\D\rfloor}|^2\\
\nonumber
&&+\sum_{i=0}^{k-1}(C^{(i+1)\D}-C^{i\D})|X_i|^2+\sum_{i=0}^{k-1}C^{(i+1)\D} m_i\\
\nonumber
&=&[-(\l_1-2\l_3-\e)\D]\sum_{i=0}^{k-1}C^{(i+1)\D}|X_i|^2\\
& &+[(1-\bar{\de})(\l_2+\l_4)\D+\e\D]\sum_{i=0}^{k-1}C^{(i+1)\D}|X_{i-\lfloor\delta(i\D)/\D\rfloor}|^2
+\sum_{i=0}^{k-1}C^{(i+1)\D} m_i,
\end{eqnarray}
which $\sum_{i=0}^{k-1}C^{(i+1)\D} m_i$ also is martingale.
Noting that
$\lfloor\delta(k\D)/\D\rfloor\leq \frac{\delta(k\D)}{\D}\leq \frac{\tau}{\D}=m$ and $C>1$.
Therefore,
\begin{eqnarray}\label{10}
\nonumber
& &\sum_{i=0}^{k-1}C^{(i+1)\D}|X_{i-\lfloor\delta(i\D)/\D\rfloor}|^2\\
\nonumber
&=&\sum_{i=-\lfloor\delta(i\D)/\D\rfloor}^{-1}C^{(i+\lfloor\delta(i\D)/\D\rfloor+1)\D}|X_i|^2
\nonumber
+\sum_{i=0}^{k-1}C^{(i+\lfloor\delta(i\D)/\D\rfloor+1)\D}|X_i|^2\\
\nonumber
& &-\sum_{i=k-\lfloor\delta(i\D)/\D\rfloor}^{k-1}C^{(i+\lfloor\delta(i\D)/\D\rfloor+1)\D}|X_i|^2\\
\nonumber
&\leq&\sum_{i=-m}^{-1}C^{(i+m+1)\D}|X_i|^2
\nonumber
+\sum_{i=0}^{k-1}C^{(i+m+1)\D}|X_i|^2\\
& &-\sum_{i=k-\lfloor\delta(i\D)/\D\rfloor}^{k-1}C^{(i+\lfloor\delta(i\D)/\D\rfloor+1)\D}|X_i|^2.
\end{eqnarray}
Substituting (\ref{10}) into (\ref{11}), we get
\begin{equation}\label{13}
C^{k\D}|X_k|^2+[(1-\bar{\de})(\l_2+\l_4)\D+\e\D]\sum_{i=k-\lfloor\delta(i\D)/\D\rfloor}^{k-1}C^{(i+\lfloor\delta(i\D)/\D\rfloor+1)\D}|X_i|^2\leq Y_k,
\end{equation}
where
\begin{eqnarray*}
Y_k
&=&|X_0|^2
+[(1-\bar{\de})(\l_2+\l_4)\D+\e\D]\sum_{i=-m}^{-1} C^{(i+m+1)\D}|X_i|^2\\
& &+\Big[-(\l_1-2\l_3-\e)\D+(1-C^{-\D})+[(1-\bar{\de})(\l_2+\l_4)\D+\e\D] C^{m\D}\Big]\sum_{i=0}^{k-1}C^{(i+1)\D}|X_i|^2\\
& &+\sum_{i=0}^{k-1}C^{(i+1)\D} m_i.
\end{eqnarray*}
Let us now introduce the function
\begin{equation}\label{3.15}
J(C,\D)=[(1-\bar{\de})(\l_2+\l_4)\D+\e\D] C^{(m+1)\D}+(1-(\l_1-2\l_3-\e)\D) C^{\D}-1.
\end{equation}
Choose $\D_1^*>0$ such that for any $\D<\D_1^*$,
$1-(\l_1-2\l_3-\e)\D>0$,
we therefore have for any $C>1$,
\begin{eqnarray*}
\frac{d}{dC}J(C,\D)=(m+1)\D[(1-\bar{\de})(\l_2+\l_4)\D+\e\D] C^{(m+1)\D-1}+(1-(\l_1-2\l_3-\e)\D)\D C^{\D-1}>0.
\end{eqnarray*}
Clearly
$$J(1)=[(1-\bar{\de})(\l_2+\l_4)+\e-\l_1+2\l_3+\e]\D,$$
since (\ref{12}), we have $J(1)<0$, which implies that there exists a unique $C^{*}_{\D}>1$ such that
$J(C_{\D}^*,\D)=0$.
We choosing $C=C_\D^*$, we therefore have
$$
Y_k=|X_0|^2+[(\l_2+2\l_3+\e)\D]\sum_{i=-m}^{-1}C_\D^{*(i+m+1)\D}|X_i|^2+\sum_{i=0}^{k-1}C_\D^{*(i+1)\D} m_i.
$$
Noting that the initial sequence $X_i<\infty$ for all $i=-m,\ldots,0$, by the semimartingale convergence Lemma,
for $C=C_{\D}^{*}$, we have
$$
\lim_{k\rightarrow\infty}Y_k<\infty \quad a.s.
$$
By (\ref{13}), we therefore have
\begin{eqnarray}\label{3.16}
\nonumber
& & \limsup_{k\rightarrow\infty} C_{\D}^{* k\D}|X_k|^2\\
\nonumber
&\leq &\limsup_{k\rightarrow\infty}[C_{\D}^{* k\D}|X_k|^2+(\l_2+2\l_3+\e)\D\sum_{i=k-\lfloor\delta(i\D)/\D\rfloor}^{k-1}C_\D^{*(i+\lfloor\delta(i\D)/\D\rfloor+1)\D}|X_i|^2]\\
&\leq &\limsup_{k\rightarrow\infty}Y_k<\infty \quad a.s.
\end{eqnarray}
Noting that $m\D=\T$, by (\ref{3.15})
\begin{equation}\label{3.17}
J(C_\D^*,\D)=[(1-\bar{\de})(\l_2+\l_4)+\e]C_{\D}^{* \T}+ \frac{1}{\D}(1-C_{\D}^{* -\D})-(\l_1-2\l_3-\e)=0.
\end{equation}
Choosing that constant $\sigma$ such that $C=e^\s$ and hence $1-C_\D=1-e^{\s\D}$.
Define
$$
\bar{J}_{\D}(\s)=[(1-\bar{\de})(\l_2+\l_4)+\e]e^{\s\T}+ \frac{1}{\D}(1-e^{-\s\T})-(\l_1-2\l_3-\e).
$$
By (\ref{3.17}) for any $\D<\D_1^*$, we have
\begin{equation}\label{3.18}
\bar{J}_{\D}(\s_{\D}^*)=0.
\end{equation}
Noting that $\lim_{\D\rightarrow 0}({1-e^{-\s\D}})/{\D}=\s$,
we have
\begin{equation}\label{3.19}
\lim_{\D\rightarrow 0}\bar{J}_{\D}(\s)=[(1-\bar{\de})(\l_2+\l_4)+\e]e^{\s\T}+ \s-(\l_1-2\l_3-\e).
\end{equation}
By the definition (\ref{14}) of $\gamma$, (\ref{3.18}) and (\ref{3.19}) yields
$$\lim_{\D\rightarrow 0}\s_{\D}^{*}=\gamma,$$
which implies that for any positive $\bar{\e}\in(0,\gamma/2)$, there exist a $\D_2^*>0$
such that for any $\D<\D_2^*$, we have
$$\s_{\D}^*>\gamma-2\e.$$
Note that (\ref{3.16}), together with the definition of $\mu_{\D}^{*}$ show that
$$
\limsup_{k\rightarrow \infty} e^{\s_{\D}^* k\D}|X_k|^2<\infty \quad a.s.
$$
So there exists a finite random variable $\g$ such that
$$
\limsup_{k\rightarrow \infty} e^{\s_{\D}^* k\D}|X_k|^2<\g\quad a.s.
$$
Beside for any $\D<\D_1^{*}\we \D_2^*$, we have
$$
0=\limsup_{k\rightarrow \infty}\frac{\log \g}{k\D}
\geq \limsup_{k\rightarrow \infty}\frac{\log(e^{\s_{\D}^{*}k \D}|X_k|)}{k\D}
=\mu_{\D}^{*}+\limsup_{k\rightarrow \infty}\frac{2\log|X_k|}{k\D},
$$
which implies that
$$
\limsup_{k\rightarrow \infty}\frac{2\log|X_k|}{k\D}\leq -\s_{\D}^*\leq -\g +2\bar{\e}\quad a.s,
$$
that is
$$
\limsup_{k\rightarrow \infty}\frac{\log|X_k|}{k\D}\leq -\s_{\D}^*\leq -\frac{\g}{2} +\bar{\e}\quad a.s
$$
as required.
The proof is hence complete.
\hfill $\Box$

We demonstrate the process of implementing the partially truncated EM by the following example.
\begin{expl}
Consider a nonlinear scalar hybrid SDDE
\begin{equation}\label{example2}
dx(t)=f(x(t),x(t-\de(t)),r(t))dt+ g(x(t),x(t-\de(t)),r(t))dB(t)
\end{equation}
Here $B(t),\de(t)$ and the Markovian chain as the same as the example \ref{example1},
where,
\begin{equation*}
f(x,y,i)=
\begin{cases}
-6x-x^5+y &\quad \text{if} \quad i=1\\
-6x-x^5+\frac{y}{1+y^2}& \quad \text{if} \quad i=2
\end{cases}
\quad and \quad
g(x,y,i)=
\begin{cases}
x^2 &\quad \text{if} \quad i=1\\
x\sin ^2y& \quad \text{if} \quad i=2
\end{cases}
\end{equation*}
\end{expl}
It can be seen that
\begin{equation*}
F_1(x,y,i)=
\begin{cases}
-6x+y &\quad \text{if} \quad i=1\\
-6x& \quad \text{if} \quad i=2
\end{cases}
\quad and \quad
G_1(x,y,i)=
\begin{cases}
0 &\quad \text{if} \quad i=1\\
0 & \quad \text{if} \quad i=2
\end{cases}
\end{equation*}
\begin{equation*}
F(x,y,i)=
\begin{cases}
-x^5 &\quad \text{if} \quad i=1\\
-x^5+\frac{y}{1+y^2}& \quad \text{if} \quad i=2
\end{cases}
\quad and \quad
G(x,y,i)=
\begin{cases}
x^2 &\quad \text{if} \quad i=1\\
x\sin^2y& \quad \text{if} \quad i=2.
\end{cases}
\end{equation*}
Choosing $\o=\8$, we then have
\begin{eqnarray*}
 2x^{T}F_1(x,y,1)+|G_1(x,y,1)|^2
=2x(-6x+y)
=-12x^2+2xy
\leq -11x^2+y^2
\end{eqnarray*}
and
\begin{eqnarray*}
 2x^{T}F_1(x,y,2)+|G_1(x,y,2)|^2 =2x(-6x)=-12x^2.
\end{eqnarray*}
Therefore, for any $i\in \mathbb{S}$, we have
\begin{equation*}
 2x^{T}F_1(x,y,i)+|G_1(x,y,i)|^2 \leq -11x^2+y^2,
\end{equation*}
that is, the (\ref{4.1}) is satisfied with $\l_1=11,\l_2=1$.
Moreover,
\begin{eqnarray*}
2x^{T}F(x,y,1)+|G(x,y,1)|^2
=2x(-x^5)+|x^2|^2
=-2x^6+x^4,
\end{eqnarray*}
but
\begin{eqnarray*}
-2x^6+x^4=-(2x^6-x^4+\frac{1}{8}x^2)+\frac{1}{8}x^2=-2x^2(x^2-\frac{1}{4})^2+\frac{1}{8}x^2 \leq \frac{1}{8}x^2.
\end{eqnarray*}
Hence
\begin{eqnarray*}
2x^{T}F(x,y,1)+|G(x,y,1)|^2\leq \frac{1}{8}x^2
\end{eqnarray*}
and
\begin{eqnarray*}
& &2x^{T}F(x,y,2)+|G(x,y,2)|^2\\
&=&2x(-x^5+\frac{y}{1+y^2})+|x\sin^2y|^2\\
&\leq& -2x^6+\frac{2xy}{1+y^2}+x^2
=-2x^6+\frac{x^2}{1+y^2}+\frac{y^2}{1+y^2}+x^2\\
&\leq& -2x^6+x^2+y^2+x^2
\leq 2x^2+y^2.
\end{eqnarray*}
Then for any $i\in \mathbb{S}$ we have
\begin{equation}
2x^{T}F(x,y,i)+|G(x,y,i)|^2 \leq 2x^2+y^2.
\end{equation}
In other words, the (\ref{4.2}) holds with the $\l_3=2,\l_4=1$. By the Theorem \ref{T4.4}, the SDDEs (\ref{example2}) is almost sure exponentially stable.
The left one of Figure 2 displays displays the almost surely asymptotic stable behaviour of the numerical solutions for the equations (\ref{example2}).
The right one of Figure 2 illustrate the almost sure exponentially stability of the numerical solution produced by the partially truncated EM method, 200 sample paths are generated.
\begin{figure}
	\centering
\includegraphics[width=6cm,height=5cm]{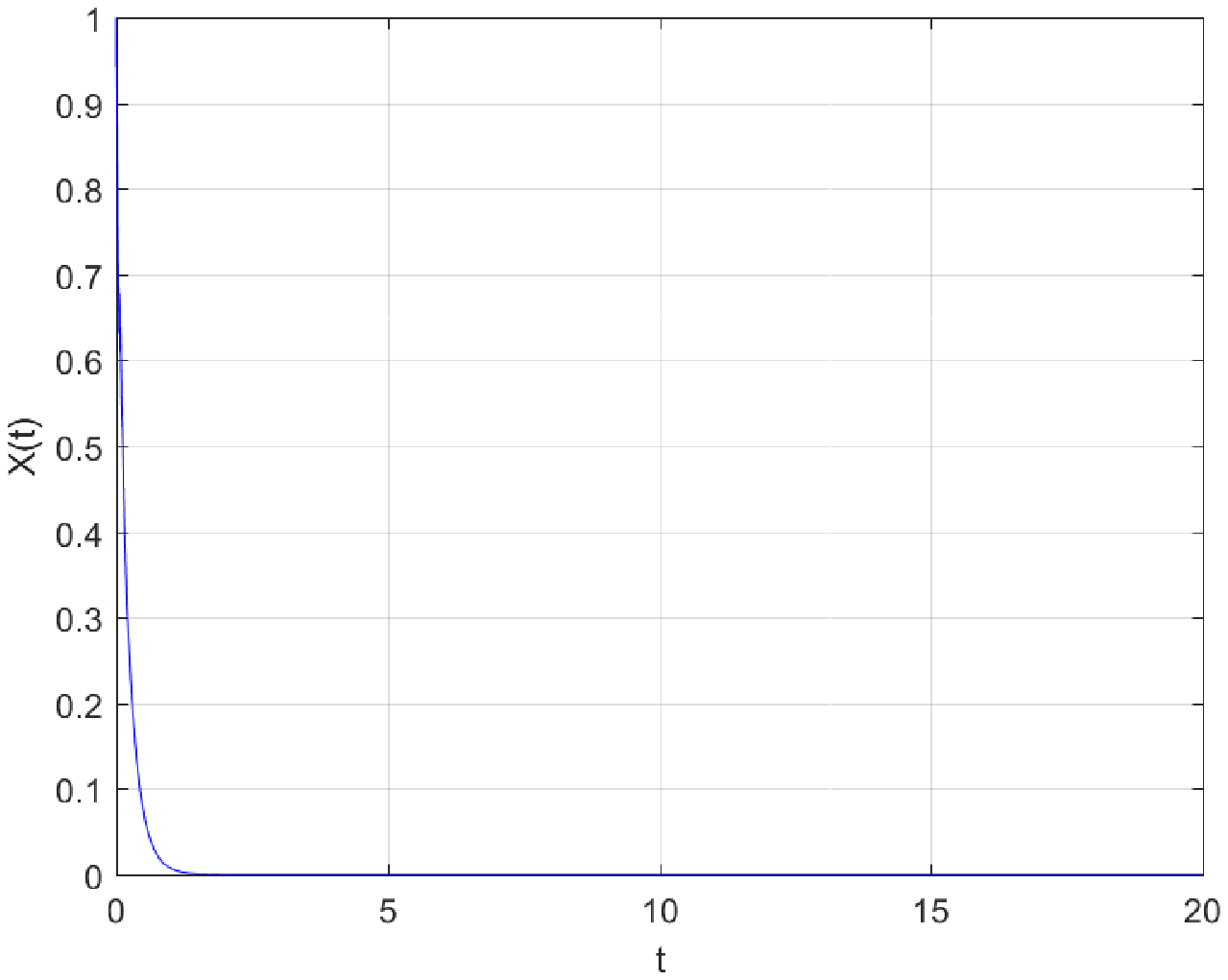}
	\includegraphics[width=6cm,height=5cm]{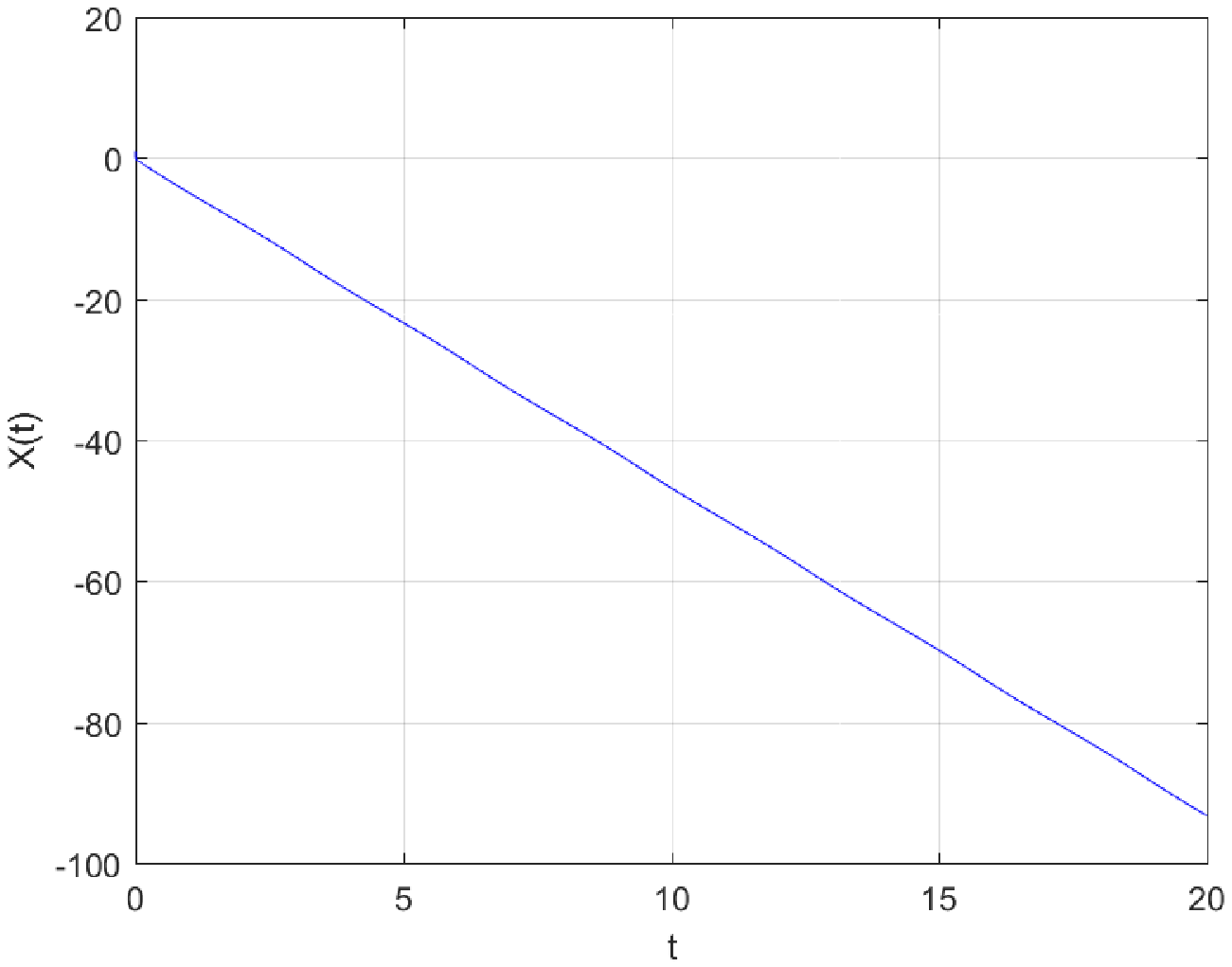}
	\caption{stability}
\end{figure}

\section*{Acknowledgements}
We thank Dr. Wei Liu and Professor Xuerong Mao for their advice.
The research of the first author is supported by the National Natural Science Foundation of China (Grant No. 11471217). 
The third author is supported by the National Natural Science Foundation of China (Grant No.  11871343).


\begin{thebibliography}{99}

\bibitem{Arn}
L.~Arnold,
 \newblock \emph{Stochastic Differential Equations: Theory and
 Applications},  John Wiley and Sons, 1972.

\bibitem{BB}
C.T.H ~Baker, and E.~Buckwar,
Numerical analysis of explicit one-step methods for stochastic delay differential equations.
\newblock\emph{LMS J. Comput. Math.}, 3(2000):315--335.

\bibitem{E.B2000}
E.~Buckwar,
Introduction to the numerical analysis of stochastic delay differnetial equations.
\newblock \emph{J. Comput. Appl. Math.}, 125(2000):297--307.

\bibitem{9}
J.~Bao, and Z.~Hou,
An analytic approximation of solutions of stochastic differential delay equations with Markovian switching.
\newblock \emph{Math. Comput. Modelling.}, 50(2009):1379--1384.


\bibitem{Guo2017}
Q.~Guo, X.~Mao, and R.~Yue,
The truncated Euler-Maruyama method for stochastic differential delay equations.
\newblock \emph{Numer Algor.}, 78(2018):599--624.

\bibitem{GLMY2017}
Q.~Guo, W.~Liu, X.~Mao, and R.~Yue,
The paryially truncated Euler-Maruyama method and its stability and boundedness.
\newblock \emph{Appl. Numer. Math.}, 115(2017):235--251.

\bibitem{12}
Q.~Guo, X.~Mao, and R.~Yue,
Almost sure exponential stability of stochastic differential delay equations.
\newblock \emph{SIAM J. Control Optim.}, 54(2016):1919--1933.


\bibitem{20}
M.~Liu, W.~Cao, and Z.~Fan,
Convergence and stability of the semi-implicit Euler method for a linear stochastic differential delay equation.
\newblock \emph{J. Comput. Appl. Math.}, 170(2004):255--268.

\bibitem{19}
R.~Li, and Y.~Hou,
Convergence and stability of numerical solutions to SDDEs with Markovian switching.
\newblock \emph{Appl. Math. Comput.}, 175(2006):1080--1091.

\bibitem{Li2010}
X.~Li, X.~Mao, and Y.~Shen,
Approximation solutions of stochastic differential delay equation with Markovian switching.
\newblock \emph{J. Differ. Equ. Appl.}, 16(2010):195--207.

\bibitem{7}
 X.~Mao, A.~Matasov, and A. B.~Piunovskiy,
Stochastic differential delay equations with Markovian switching.
\newblock \emph{Bernoulli.}, 6(2000):73--90.

\bibitem{MR05}
X.~Mao, and M.J.~Rassias,
Khasminskii-type theorems for stochastic differential delay equations.
\newblock \emph{J. Sto. Anal. Appl.}, 23(2005):1045--1069.

\bibitem{MY06}
X.~Mao, and C.~Yuan,
\newblock \emph{Stochastic Differential Equations with Markovian Switching},
Imperial College Press, 2006.

\bibitem{MYZ}
X.~Mao, C.~Yuan, and J.~Zou,
Stochastic differential delay equations of population dynamics.
\newblock \emph{J. Math. Anal. Appl.}, 304(2005):296--320.

\bibitem{MS}
X.~Mao, and S.~Sabanis,
Numerical solutions of stochastic differential delay equations under local Lipschitz condition.
\newblock \emph{J. Comput. Appl. Math.}, 151(2003):215--227.

\bibitem{mao2015truncated}
X.~Mao,
The truncated euler--maruyama method for stochastic differential equations.
\newblock \emph{J. Comput. Appl. Math.}, 290(2015):370--384.

\bibitem{mao2016convergence}
X.~Mao,
Convergence rates of the truncated euler--maruyama method for stochastic differential equations.
\newblock \emph{J. Comput. Appl. Math.}, 296(2016):362--375.

\bibitem{M2002}
X.~Mao,
A note on the LaSalle-type theorems for stochastic differential delay equations.
\newblock \emph{J. Math. Anal. Appl.}, 268(2002):125--142.

\bibitem{M11}
X.~Mao,
Numerical solutions of stochastic differential delay equations under the generalized Khasminskii-type conditions.
\newblock \emph{Appl. Math. Comput.}, 217(2011):5512--5524.

\bibitem{mao1999}
X.~Mao,
Stability of stochastic differential equations with Markovian switching.
\newblock \emph{Stoch. Proc. Appl.}, 79(1999):45--67.

\bibitem{mao2002}
X.~Mao,
Exponential stability of stochastic delay interval systems with Markovian switching.
\newblock \emph{IEEE Trans. Automat. Control.}, 47(2002):1604--1612.

\bibitem{M02}
X.~Mao,
A note on the LaSalle-type theorems for stochastic differential delay equations.
\newblock \emph{J. Math. Anal. Appl.}, 268(2002):125--142.

\bibitem{Mao2007}
X.~Mao,
Exponential stability of equidistant Euler-Maruyama approximations of stochastic differential delay equations.
\newblock \emph{J. Comput. Appl. Math.}, 200(2007):297--316.

\bibitem{1}
 X.~Mao,
Stability of stochastic differential equations with Markovian switching.
\newblock \emph{Stoch. Proc. Appl.}, 79(1999):45--67.

\bibitem{21}
M.~Song, L.~Hu, X.~Mao, and L.~Zhang,
Khasminskii-type theorems for stochastic functional differential equations.
\newblock \emph{Discrete Contin. Dyn. Syst. Ser. B}, 18(2013):1697--1714.

\bibitem{8}
Y. Shen, M.~Jiang, and X.~Liao,
The LaSalle-type theorem for stochastic differential delay equations with Markovian switching.
\newblock \emph{Dyn. Contin. Discrete Impuls. Syst. Ser. A Math. Anal.}, 13(2006):1254--1262.

\bibitem{wu2010}
F.~Wu, X.~Mao, and L.~Szpruch,
Almost sure exponential stability of numerical solutions for stochastic delay differential equations.
\newblock \emph{Number. Math.}, 115(2010):681--697.

\bibitem{5}
L.~Wang, and F.~Wu,
Existence, uniqueness and asymptotic properties of a class of nonlinear stochastic differential delay equations with Markovian switching.
\newblock \emph{Stoch. Dynam.}, 9(2009):253--275.

\bibitem{15}
C.~Yuan, and W.~Glover,
Approximate solutions of stochastic differential delay equations with Markovian switching.
\newblock \emph{J. Comput. Appl. Math.}, 194(2006):207--226.

\bibitem{2}
F.~Zhu, Z.~Han, and J.~Zhang,
Stability analysis of stochastic differential equations with Markovian switching.
\newblock \emph{Syst. Control. Lett.}, 61(2012):1209--1214.

\bibitem{22}
W.~Zhang, M.~Song, and M.~Liu,
Strong convergence of the partially truncated Euler-Maruyama method for a class of stochastic differential delay equations.
\newblock \emph{J. Comput. Appl. Math.}, 335(2018):114--128.

\bibitem{3}
E.~Zhu, X.~Tian, and Y.~Wang,
On $p$-th moment exponential stability of stochastic differential equations with Markovian switching and time-varying delay.
\newblock \emph{ J. Inequal. Appl.}, 1(2015):1--11.

\bibitem{18}
S.~Zhou,
Strong convergence and stability of backward Euler--Maruyama scheme for highly nonlinear hybrid stochastic differential delay equation.
\newblock \emph{Calcolo: A quarterly on Numerical Analysis and Theory of Computation.}, 52(2015):445--473.

\end{thebibliography}
\end{document}